\documentclass[fleqn,a4paper,10pt]{amsart}

\usepackage{graphicx}
\usepackage{tikz}
\usetikzlibrary{shapes,decorations.pathmorphing}
\usepackage{amsmath}
\usepackage{amssymb,amsthm}
\usepackage{comment}
\usepackage{hyperref}

\def\cB{\mathcal{B}}
\def\cE{\mathcal{E}}

\def\cM{\mathcal{M}}

\theoremstyle{definition}

\def\tb{t_{\bullet}}
\def\tw{t_{\circ}}

\def\cBnd{\cB_n^{(d)}}
\def\cM{\mathcal{M}}
\def\cMnd{\cM_n^{(d)}}
\def\cBhnd{\widehat{\cB}_n^{(d)}}
\def\cMhnd{\widehat{\cM}_n^{(d)}}
\def\tb{t_{\bullet}}
\def\hF{\widehat{F}}
\def\hJ{\widehat{J}}
\def\hQ{\widehat{Q}}
\def\tw{t_{\circ}}
\def\wF{\widehat{F}}
\def\hM{\widehat{M}}

\title[Comparing ensembles of quadrangulations via
continued fractions]{Comparing two statistical ensembles of quadrangulations:
a continued fraction approach}
\author{\'Eric Fusy}
\address{LIX, \'Ecole Polytechnique, 91120 Palaiseau, France. Academic year 2014-2015: PIMS-CNRS, University of British Columbia, Vancouver, BC, Canada.}
\email{fusy@lix.polytechnique.fr}

\author{Emmanuel Guitter}
\address{Institut de Physique Th\'eorique, CEA, IPhT, 91191 Gif-sur-Yvette, France, CNRS, URA 2306}
\email{emmanuel.guitter@cea.fr}

\begin{document}
\maketitle

\begin{abstract}
We use a continued fraction approach to compare two statistical ensembles of quadrangulations with a boundary,
both controlled by two parameters. In the first ensemble, the quadrangulations are bicolored and the parameters control their numbers of
vertices of both colors. In the second ensemble, the parameters control instead the number of vertices which are local maxima 
for the distance to a given vertex, and the number of those which are not. Both ensembles may be described either by their (bivariate) 
\emph{generating functions at fixed boundary length} or, after some standard slice decomposition, by their (bivariate) 
\emph{slice generating functions}. 
We first show that the fixed boundary length generating functions are in fact equal for the two ensembles. We then show that
the slice generating functions, although different for the two ensembles, simply correspond to two different ways of encoding 
the same quantity as a continued fraction. This property is used to obtain explicit expressions for the slice generating functions in a constructive way.
\end{abstract}

\section{Introduction}
\label{sec:introduction}
The study of planar maps has given rise in the recent years to a lot of remarkable enumeration results. A
particularly fruitful approach consists in taking advantage of bijections between maps and tree-like objects called mobiles. This technique,
initiated by Schaeffer \cite{SchPhD,CMS09} (reinterpreting a bijection by Cori and Vauquelin \cite{CoriVa}) was extended
in many different directions \cite{BDG04,AmBudd} to deal with various refined map enumeration problems. Besides mobiles, 
another, slightly different, view on the problem consists in decomposing the maps into so-called \emph{slices}, which are particular pieces of maps
with nice combinatorial properties \cite{BG12}. In particular, the generating functions for these slices were shown to obey discrete 
integrable systems of equations and in most cases, a solution of these equations could be obtained explicitly. Moreover, the slice
decomposition of a map is intimately linked to its geodesic paths and the knowledge of slice generating functions
directly gives explicit answers to a number of questions regarding the statistics of distances between 
random points within maps \cite{GEOD,BDG04,BFG,FG14}.    

A particularly important discovery was made in \cite{BG12} where it was shown that slice generating functions happen to be
simple coefficients in a suitable \emph{continued fraction expansion} of standard map generating functions,
making {\it de facto} a connection between the distance statistics within maps and some more global properties.
On a computational point of view, this discovery provided a constructive way to obtain explicit solutions  
for the integrable systems at hand, by taking advantage of known results on continued fractions.

Quite recently, the slice decomposition technique was used in \cite{FG14} to describe the distance statistics of general families of bicolored maps, 
and, in particular, of \emph{bicolored quadrangulations}, with some simultaneous control on the numbers of vertices of both colors.
Explicit expressions for the corresponding bivariate slice generating functions were obtained in a constructive way,
leading in particular to explicit formulas for the distance dependent two-point function within bicolored quadrangulations.
Remarkably, the expressions found for slice generating functions are very similar to those
obtained (via a mobile formalism) in another problem of quadrangulations considered in \cite{AmBudd}. There, the discrimination
between vertices no longer relies on their color but rather on their status with respect to the graph distance from a fixed origin vertex. 
Vertices namely come in two types: those, called \emph{local maxima} which are further from the origin than all their neighbors, and the others. Bivariate
slice generating functions can be defined so as to keep some independent control on the numbers of both types of vertices
after the slice decomposition. Explicit expressions for these new bivariate slice generating functions were then \emph{guessed}
in \cite{AmBudd} and, as just mentioned, their structure is very similar to that of their bicolored counterparts.

The aim of this paper is twofold: first, we establish a strong connection between the problem of quadrangulations with a control on the vertex 
color, as discussed in \cite{FG14}, and that of quadrangulations with a control on local maxima, as discussed in \cite{AmBudd}. Then, we use a continued fraction formalism to re-derive, now in a constructive way,
the explicit expressions found in \cite{AmBudd}.

The paper is organized as follows: in Sect.\ref{sec:bivariate}, we introduce the two ensembles of quadrangulations 
that we want to compare and define their generating functions at fixed boundary length.
We then derive our first fundamental result which states that the fixed boundary length generating functions are in fact equal 
for the two ensembles. Sect.~\ref{sec:contfrac} presents the slice decomposition of the quadrangulations at hand and 
shows that the corresponding slice generating functions may be obtained as coefficients of the same quantity, once expanded as a continued fraction
in two different ways.
In Sect.~\ref{sec:recureq}, we recall the integrable systems which determine the slice generating functions of both ensembles as well as the
explicit solutions of these systems obtained in \cite{FG14} and \cite{AmBudd}. Sect.~\ref{sec:extraction} deals with results on continued fractions,
and shows in particular how to extract their coefficients from those obtained via a standard series expansion. In one of the ensembles
that we consider, the knowledge of this series expansion is not sufficient to get all the slice generating functions,
a process which requires \emph{the knowledge of some additional quantity}. 
An explicit expression for this latter quantity is conjectured in Sect.~\ref{sec:recovering}, based on simplifications observed
in the case of finite continued fractions, and we then show how it allows to recover the explicit
formulas for the slice generating functions found in \cite{AmBudd}. 
Sect.~\ref{sec:conserved} deals with another aspect of our problem, the existence of invariants, the so-called \emph{conserved quantities},
as expected for discrete integrable systems. We show how to derive these invariants combinatorially for both 
ensembles and again emphasize the deep similarity existing between the conserved quantities for the two ensembles.
We gather our concluding remarks in Sect.~\ref{sec:conclusion}. Some side results or technical derivations are presented in Appendices A and B.

\section{An equality between two bivariate generating functions for quadrangulations with a boundary}
\label{sec:bivariate}
The aim of this section is to compare the generating functions of planar \emph{quadrangulations with a boundary} 
weighted in two different ways, each of these weighting being \emph{bivariate}, i.e.\ involving two independent parameters.
As we shall see below, these two weightings, although fundamentally different, are intimately linked and some of the
associated generating functions turn out to be equal. 

\subsection{Two bivariate generating functions for quadrangulations with a boundary}
\label{sec:weightings}
Recall that a planar quadrangulation with a boundary denotes a connected graph embedded on the sphere 
which is \emph{rooted}, i.e.\  has a marked oriented edge (the root-edge) and is such that all 
its inner faces, i.e.\ all the faces except that lying on the right of the root-edge, have degree $4$. As for the external face, which is the face 
lying on the right of the root-edge, its degree is arbitrary (but necessarily even). As customary, the origin of the 
root-edge will be called the root-vertex. Let us now consider two particular different ways to assign weights to 
these maps.
\begin{figure}
\begin{center}
\includegraphics[width=12cm]{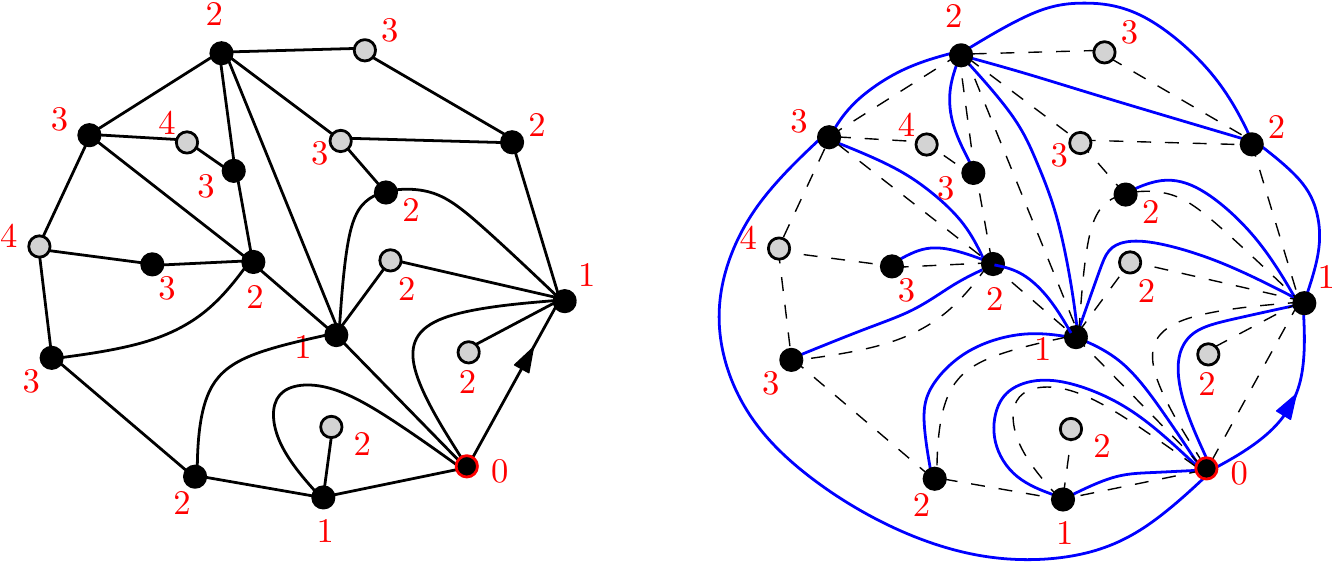}
\end{center}
\caption{Left: an example of rooted quadrangulation with a boundary of length $10$, where each vertex
is labelled by its distance to the root-vertex. The local maxima for this distance are indicated in gray.
Right: the associated rooted general map with a bridgeless boundary of length $5$, obtained by applying the Ambj\o rn-Budd rule
within each inner face (i.e.\ connecting the two corners within the face followed clockwise by a larger label) and,
in the external face, connecting \emph{cyclically} those corners followed by a larger label counterclockwise around the map.}
\label{fig:quadquad1}
\end{figure}
\vskip .5cm
\noindent \checkmark \emph {First weighting: bicoloring the map.}
Since all their faces have even degree, planar quadrangulations with a boundary may be naturally bicolored
in black and white in a unique way, by assigning the black color to the root-vertex and demanding that no two adjacent vertices have the same color.
We way then enumerate these quadrangulations by assigning a weight $t_\bullet$ to each black vertex and a weight $t_\circ$ to each white vertex.
For convenience, the root-vertex receives a weight $1$ instead of $t_\bullet$. We shall then denote by $F_n\equiv F_n(t_\bullet,t_\circ)$ 
the corresponding generating function for these maps with a \emph{boundary length} $2n$, i.e.\ with an external face of degree $2n$.
\vskip .5cm 
\noindent \checkmark \emph {Second weighting: distinguishing local maxima of the distance}. Our second weighting consists in giving a special role
to the local maxima of the distance from the root-vertex. More precisely, we may label each vertex $v$ of the quadrangulation by its 
graph distance $d(v)$ from the root-vertex and look for the \emph{local maxima} of this labeling, i.e.\ those vertices $v$ having only neighbors
with label $d(v)-1$ (note that in all generality, neighbors of a vertex $v$ may only be at distance $d(v)-1$ or $d(v)+1$ from the root-vertex).
We decide to give a weight $t_\circ$ to local maxima and a weight $t_\bullet$ to the other vertices (see Fig.~\ref{fig:quadquad1}-- left). As before, the root-vertex
receives a weight $1$ instead of $t_\bullet$ (note that the root-vertex can never be a local maximum). We shall call 
$J_n\equiv J_n(t_\bullet,t_\circ)$ the generating function for these maps with a boundary of length $2n$. 

The generating functions $F_n(t_\bullet,t_\circ)$ and $J_n(t_\bullet,t_\circ)$  may be understood as formal power series
in $t_\bullet$ and $t_\circ$, giving rise to convergent series for small enough 
$t_\bullet, t_\circ$.
The first weighting is quite natural and was described in detail in \cite{FG14}. We shall recall some of the corresponding results below.
As for the second weighting, it may seem more artificial but, as explained in \cite{AmBudd}, it arises naturally in two contexts:
first, letting $t_\circ\to 0$ (i.e\ keeping the linear term in $t_\circ$) is a way to suppress local maxima of the distance, selecting quadrangulations arranged into layers between the root-vertex and a unique local maximum. These so-called Lorentzian or causal structures display a very different statistics
from that of arbitrary quadrangulations \cite{AmBudd}. As recalled below, the second weighting also arises naturally when
enumerating \emph{general planar maps} with a control on both their numbers of vertices and faces \cite{AmBudd}.  

\subsection{Equality of generating functions}
\label{sec:equality}
Let us now prove a first fundamental equality, namely that
\begin{equation}
J_n(t_\bullet,t_\circ)=F_n(t_\bullet,t_\circ)\ .
\label{eq:fundequal}
\end{equation}
To this end, let us recall the so-called Ambj\o rn-Budd bijection of \cite{AmBudd} between quadrangulations and general maps, slightly 
adapted to the case of quadrangulations with a boundary according to the rules of \cite{BFG}. Starting with our quadrangulation with a 
boundary and labeling each vertex $v$ by its distance $d(v)$ from the root-vertex, we associate to each inner face an edge as follows 
(see Fig.~\ref{fig:quadquad1}-- right): looking at the corners\footnote{Recall that a corner is an angular sector between two successive half-edges 
around a given vertex. The label of a corner is that of the incident vertex.} clockwise around the face, exactly two corners are followed
 by a corner with larger label. 
We connect these two corners by an edge lying inside the original face. As for the external face, looking again at the corner labels 
clockwise around the face, i.e.\ counterclockwise around the rest of the map, exactly $n$ corners, including the root-corner (lying immediately to 
the right of the the root-edge) are followed 
by a corner with larger label. We connect the $n$ corners of this ensemble \emph{cyclically} clockwise around the map, each edge connecting two successive corners in the ensemble (see Fig.~\ref{fig:quadquad1}). Finally, we mark and orient away from the root-vertex the edge connecting
the root-corner to its successor. As explained in \cite{AmBudd} (and its extension \cite{BFG}), the obtained edges form a rooted planar map together with those vertices of the original
quadrangulation \emph{which were not local maxima} for the distance from the root-vertex. Each inner face of this map surrounds exactly one 
of the original local maxima, which get disconnected in the construction.
More precisely, from \cite{AmBudd,BFG}, the above transformation provides a bijection between planar quadrangulations with a boundary
of length $2n$ and rooted planar general (i.e.\ with faces of arbitrary degrees) maps with a \emph{bridgeless} boundary\footnote{Strictly speaking, the extension \cite{BFG} of \cite{AmBudd} shows that corners followed by
a smaller label should be connected \emph{cyclically} within all faces, including the inner faces, so that the resulting object is a hypermap, made of alternating black and white faces with, in our case, all black faces of degree $2$ but one, of degree $n$, which we choose as external face. 
The Ambj\o rn-Budd construction that we use here is recovered by squeezing all inner black faces, of degree $2$, into simple edges 
while the external black face of degree $n$ becomes the external face of the map. As for any face of a hypermap, its boundary is then necessarily without bridge.} of length $n$ (i.e.\ with external face -- lying on the right of the root-edge -- of degree $n$ and without bridge).
The vertices of the quadrangulation which are not local maxima for the distance from the root are in one-to-one correspondence with the vertices of the general map while the vertices of the quadrangulation which are local maxima are in one-to-one correspondence with the inner faces of the general map.

We may thus interpret $J_n(t_\bullet,t_\circ)$ as the generating function for rooted planar general maps with a bridgeless boundary 
of length $n$, weighted by
$t_\bullet$ per non-root-vertex 
and $t_\circ$ per inner face.  

As for the label $d(v)$ of a vertex $v$ retained in this new map, 
it precisely corresponds to the \emph{oriented} graph distance from the root-vertex to $v$ on the new map, using paths oriented from the root-vertex to $v$ which respect 
the following edge orientation\footnote{In the underlying hypermap structure, the labels correspond to the distance using oriented paths going clockwise around the black faces. Squeezing the inner black faces of degree $2$ results in simple edges 
oriented both ways, while the boundary-edges remain oriented oneway only.}: 
all edges are oriented both ways except for the boundary-edges 
(i.e.\ the edges incident to the external face) which are oriented counterclockwise around the map.

\begin{figure}
\begin{center}
\includegraphics[width=12cm]{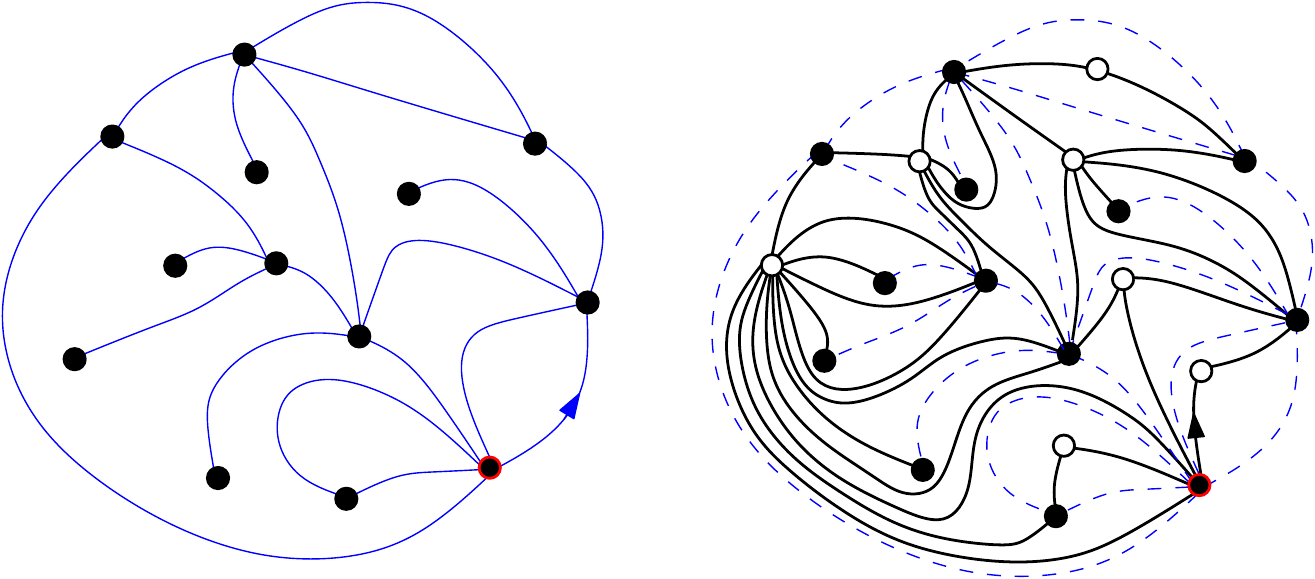}
\end{center}
\caption{Left: the rooted map on the right of Fig.~\ref{fig:quadquad1}, with a bridgeless boundary of length $5$. Right:
The associated rooted bicolored quadrangulation with a boundary of length $10$ obtained by inserting a white vertex at the center 
of each inner face of the map and connecting it to all the incident vertices around the face.}
\label{fig:quadquad2}
\end{figure}

Forgetting about distances and labels, we may now use a standard construction to rebuild a quadrangulation with a boundary out of 
our general map. Coloring the vertices of the map in black, we simply add a white vertex within each inner face
and connect it to all the corners within the face\footnote{Note that it is crucial that the boundary of the map be bridgeless
for the obtained object to be connected.}. 
By doing so, we get a bicolored quadrangulation with a boundary twice larger as that of the general map we started from, which we root by picking 
the edge leaving the root-vertex within the corner immediately to the left of the root-edge of the general map, and orienting it from its black to its white extremity (see Fig.~\ref{fig:quadquad2}-- right). Again this construction provides a bijection between rooted planar general maps with a bridgeless boundary of length $n$ and
planar quandrangulations with a boundary of length $2n$, equipped with their (unique) 
bicoloration  
as defined in the previous section. The vertices of the general map are in one-to-one correspondence with the black vertices of the quadrangulation while the inner faces of the map are in one-to-one correspondence with the white vertices of the quadrangulation.

We may thus interpret $F_n(t_\bullet,t_\circ)$ as the generating function for rooted planar general maps with a bridgeless boundary of length $n$, weighted
by $t_\bullet$ per non-root-vertex and $t_\circ$ per face. Eq.~\eqref{eq:fundequal} follows. 

 \section{Slice decomposition and continued fractions}
\label{sec:contfrac}
\subsection{Slice decomposition for maps enumerated by \boldmath{$F_n$}}
\label{sec:contfracFn}
As explained in \cite{FG14}, the quadrangulations (with a boundary) 
enumerated by $F_n$ may be decomposed into \emph{slices}
by some appropriate cutting procedure.

\begin{figure}
\begin{center}
\includegraphics[width=8cm]{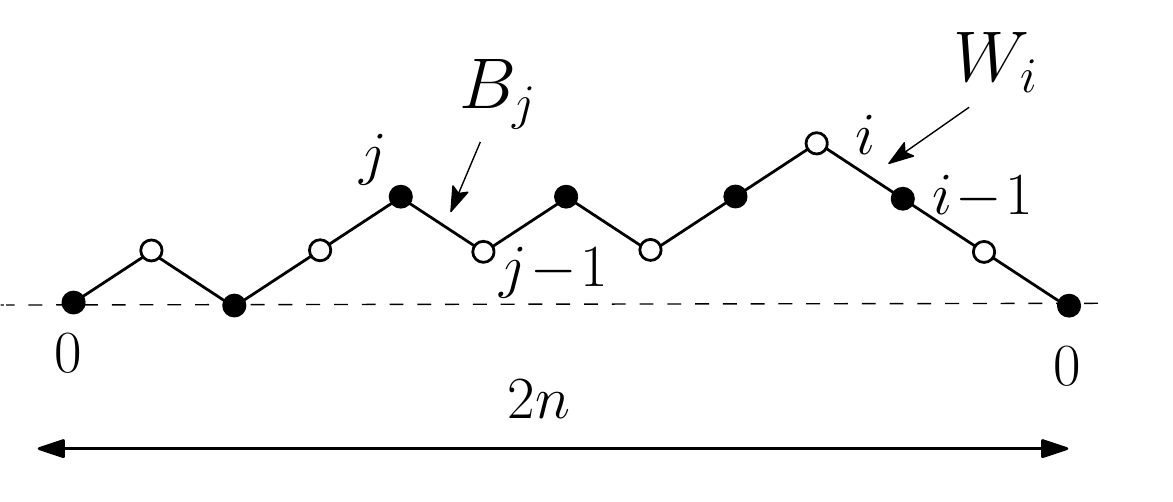}
\end{center}
\caption{An example of directed path of length $2n=12$, made of elementary steps with height difference $\pm 1$, 
starting and ending at height $0$ and remaining (weakly) above height $0$. The path is naturally colored in black and white.
To compute $F_n$, we must sum over such bicolored paths
with a weight $B_i$ (resp.\ $W_i$) assigned to each descending step $i\to i-1$  starting from a black height (resp. white height).}
\label{fig:pathBiWi}
\end{figure}

\begin{figure}
\begin{center}
\includegraphics[width=10cm]{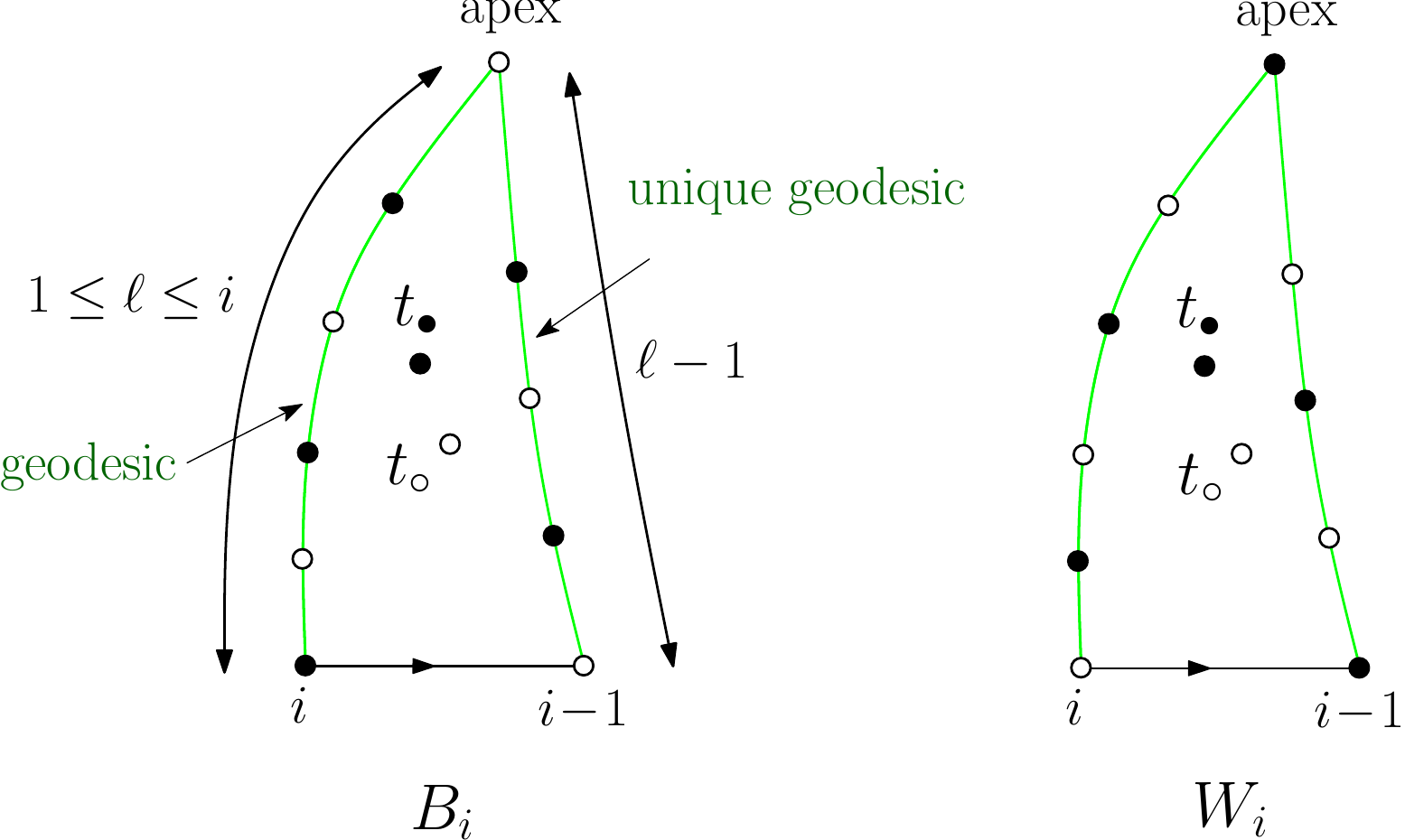}
\end{center}
\caption{A schematic picture of an $i$-slice contributing to $B_i$ (left) and to $W_i$ (right).}
\label{fig:islice}
\end{figure}

Labeling each boundary-vertex $v$ by 
$d(v)$, the sequence of corner labels, read counterclockwise around the map
starting from the root-corner, forms a directed path of length $2n$, made of elementary steps with height difference $\pm 1$, 
starting and ending at height 
$0$ and remaining (weakly) above height 
$0$ (see Fig.~\ref{fig:pathBiWi}). Drawing, for each boundary-vertex $v$, its \emph{leftmost}
geodesic (shortest) path to 
the root-vertex and cutting along these geodesics results into a decomposition of the map into pieces,
called slices. More precisely, to each descending step $i\to i-1$ of the path  corresponds an $i$-slice, defined as follows
(see \cite{FG14} for details): it is a rooted map whose boundary is made of three parts (see Fig.~\ref{fig:islice}): 
(i) its base consisting of a single root-edge, (ii) a left boundary of length $\ell$ with $1\leq \ell \leq i$ connecting the 
origin of the root-edge to another vertex, the \emph{apex} and which is a geodesic path within the slice, and (iii) 
a right boundary of length $\ell-1$ connecting the endpoint of the root-edge to the apex, and which is the unique geodesic path 
within the slice between these vertices. The left and right boundaries do not meet before reaching the apex
(which by convention is considered as part of of the right boundary only). 
As a degenerate case when $\ell=1$, the left boundary may stick to the base, in which case the slice is reduced to a single 
root-edge.

At this level, it is interesting to note that the distance $d(v)$ from the root-vertex to
any vertex $v$ in the quadrangulation
is directly related to its distance $d_s(v)$, within the $i$-slice it lies in, from the apex of this slice via
\begin{equation}
d(v)=d_s(v)+i-\ell
\label{eq:dds}
\end{equation} 
if $\ell$ is the length of the left boundary of the slice. Indeed, it is clear by construction of the slices that either  
the root-vertex is the apex of the $i-$slice at hand
or it does not belong to the slice at all and any path from $v$ to 
this root-vertex must first reach one of the boundaries of the slice (possibly at the apex). In the first case, we have  
$d(v)=d_s(v)$ and  
$i=\ell$ so that \eqref{eq:dds} holds. In the second case, since the slice boundaries are part of geodesic paths to  
the root-vertex, $d(v)$ is equal to $d_s(v)$ plus the distance from the apex of the slice to 
the root-vertex. In other words, 
$d(v)-d_s(v)$ has a constant value within the $i$-slice, which is obtained by taking for $v$ the origin of the root-edge of the slice, 
namely $d(v)-d_s(v)= i-\ell$, and \eqref{eq:dds} follows. Note that, in an $i$-slice, $i$ only acts as an upper bound on the length $\ell$ of the left boundary.
The vertices $v$ of an $i$-slice may then be labelled by non-negative integers in two natural ways: either by their distance $d_s(v)$ to the apex or by 
this distance plus $i-\ell$.

Let us call $B_i\equiv B_i(t_\bullet,t_\circ)$ (resp.\ $W_i\equiv W_i(t_\bullet,t_\circ)$) the generating function for $i$-slices whose root-vertex is black (resp.\ white), with a weight $t_\bullet$ per black vertex and $t_\circ$ per white vertex \emph{except for the vertices of the right boundary} (including
the endpoint of the root-edge and the apex) which receive a weight $1$~\footnote{The fact that these vertices receive a weight $1$
is to avoid double weighting upon re-gluing the slices into a quadrangulation. Indeed, all these vertices are 
already part of a left boundary, except for the root-vertex. In the end, only the  
root-vertex gets a weight $1$, as wanted.}.
The slice decomposition implies that \cite{FG14}
\begin{equation}
F_n=Z_{0,0}^{+}(2n;\{B_i\}_{i\geq 1},\{W_i\}_{i\geq 1})\ ,
\label{eq:Fnexpr}
\end{equation}
where $Z_{0,0}^{+}(2n;\{B_i\}_{i\geq 1},\{W_i\}_{i\geq 1}$ denotes the generating function of paths of length $2n$,
made of elementary steps with height difference $\pm 1$, colored alternatively in black and white, starting and ending at black height  
$0$ and 
remaining (weakly) above height 
$0$, with each descending step from a black height $i$ to a white height $i-1$ weighted by
$B_i$ and each descending step from a white height $i$ to a black height $i-1$ weighted by
$W_i$ (see Fig.~\ref{fig:pathBiWi}). 

The set of identities \eqref{eq:Fnexpr} for all $n>0$ can be summarized into the continued fraction expression 

\begin{equation}
F(z)\equiv \sum_{n\geq 0} F_n z^n=\frac{1}{\displaystyle{1-z \frac{W_1}{\displaystyle{1-z \frac{B_2}{\displaystyle{1-z \frac{W_3}{\displaystyle{1-z \frac{B_4}{\displaystyle{ 1- \cdots}}}}}}}}}}
\label{eq:contfrac}
\end{equation}
with the convention that $F_0=0$ and where $F(z)=F(z;t_\bullet,t_\circ)$ implicitly depends on $t_\bullet$ and $t_\circ$.
 
\subsection{Slice decomposition for maps enumerated by \boldmath{$J_n$}}
\label{sec:contfracJn}
Let us now play the same game with maps enumerated by $J_n$, which are the same maps 
as those enumerated by $F_n$, but now with the second weighting. 
We may again apply the same slice decomposition, resulting in the same $i$-slices as before.
More precisely, labeling each boundary-vertex $v$ by $d(v)$ gives rise to a path of length $2n$
(from height $0$ to height $0$, remaining above height $0$) and each descending step $i\to i-1$ 
gives rise to an $i$-slice. To assign the second weighting to the quadrangulation, we must 
label each vertex $v$ of the $i$-slice by its distance $d(v)$ (in the quadrangulation) from the root-vertex 
of the quadrangulation. As explained above, if the $i$-slice has a left boundary length $\ell$ ($1\leq \ell \leq i$),
this amounts to label $v$ by $d_s(v)+i-\ell$ where $d_s(v)$ is its distance (within the slice)
from the apex of the slice. To recover the correct weights, we must first give weight $1$ to all the vertices of the right boundary 
(including the endpoint of their root-edge and the apex) in order to avoid double weightings after regluing the slices. 
As for the vertices lying on the left boundary of the slice and \emph{different from the root-vertex of the slice}, 
they cannot, as part of a geodesic of the original quadrangulation, be local maxima as they have a
neighbor with larger label along the geodesic path. They receive a weight $t_\bullet$ accordingly.
Considering now vertices $v$ lying strictly within the slice, they have all their original neighbors lying in the slice and, from 
\eqref{eq:dds}, are local maxima for the distance $d(v)$ if and only if they are local maxima for
the distance $d_s(v)$ within the slice. For all these vertices, we may thus use the distance $d_s(v)$ within the slice 
to detect the local maxima, and give them the weight $t_\circ$, while non-local maxima for $d_s(v)$ get the weight $t_\bullet$.
The last vertex to consider is the origin $w$ of the root-edge of the slice: for this vertex to be a local maximum of the original 
distance $d(v)$, it must both be a local maximum of the distance $d_s(v)$ within the slice and have no neighbor with 
larger label after regluing. Now two situations may occur: either the boundary-vertex preceding $w$ along the boundary (oriented
counterclockwise around the quadrangulation) has label $d(w)-1$ and then the slice at hand sticks to the boundary so that
$w$ has all its neighbors within the slice. Then $w$ is a local maximum for $d(v)$ if and only if it is a local maximum for $d_s(v)$.
Or  the boundary-vertex preceding $w$ has label $d(w)+1$ is which case $w$ is not a local maximum for $d(v)$, irrespectively
of whether or not it is one for $d_s(v)$.  

\begin{figure}
\begin{center}
\includegraphics[width=10cm]{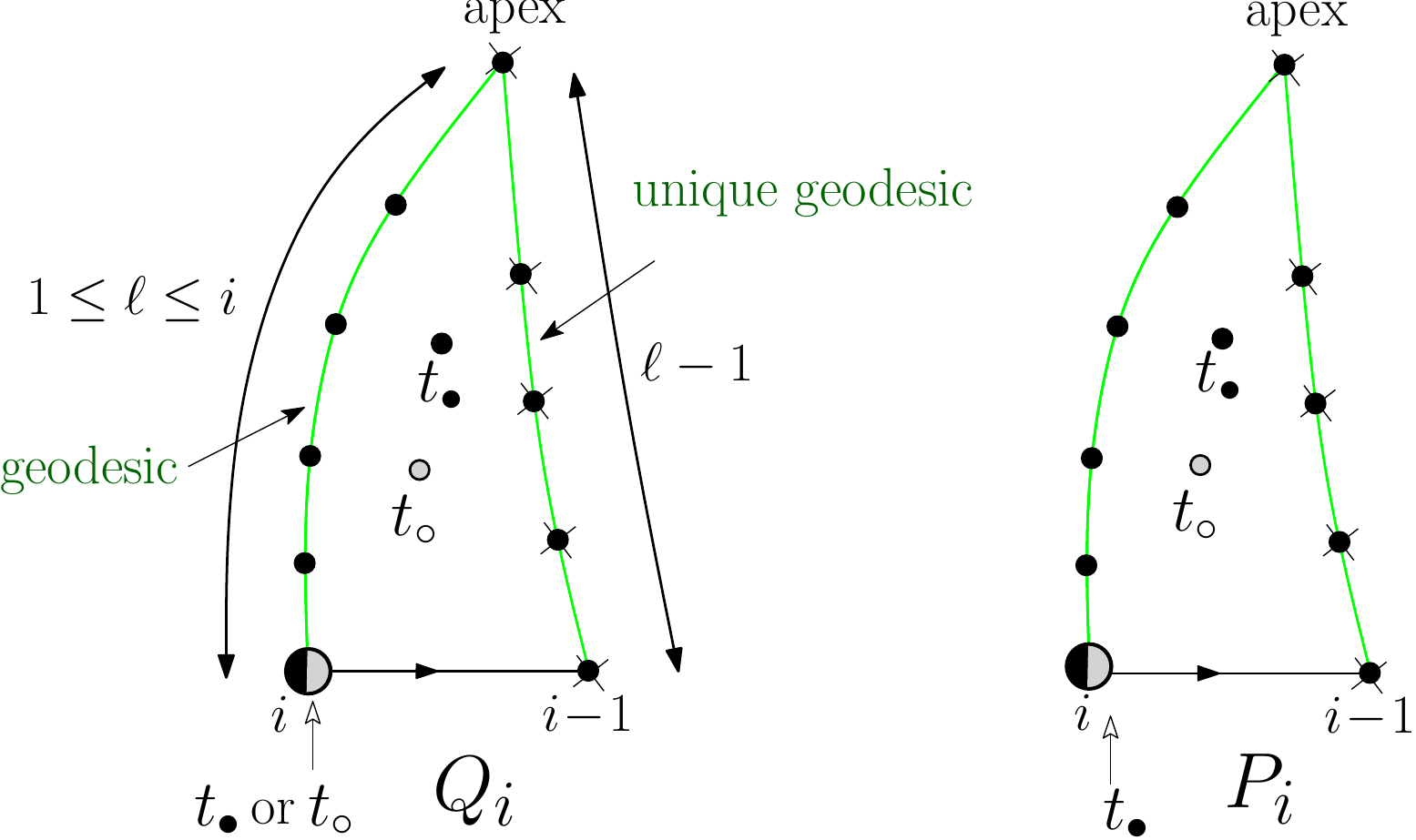}
\end{center}
\caption{A schematic picture of an $i$-slice contributing to $Q_i$ (left) and to $P_i$ (right). Local maxima of the distance from the apex 
are indicated in gray and non local maxima in black. The root vertex
may be a local maximum or not, and receives the weight $t_\circ$ or $t_\bullet$ accordingly in $Q_i$, while it always gets
the weight $t_\bullet$ in $P_i$.}
\label{fig:newislice}
\end{figure}

To summarize, we are led to consider two different generating functions for $i$-slices. In the first generating function 
$Q_i\equiv Q_i(t_\bullet,t_\circ)$, all the vertices of the $i$-slice receive a weight $t_\circ$ or $t_\bullet$ according to whether or not
they are a local maximum \emph{for the distance $d_s(v)$ from the apex within the slice}
(in particular the vertices of the left boundary different from the root-vertex of the slice get a weight $t_\bullet$ as wanted), except for the vertices of the right boundary 
(including the endpoint of their root-edge and the apex)  which receive a weight $1$. In the second generating function
$P_i\equiv P_i(t_\bullet,t_\circ)$, we assign exactly the same weights as in $Q_i$, except for the root-vertex of the
slice, which gets the weight $t_\bullet$ irrespectively of whether or not it is a local maximum for $d_s(v)$ (see Fig.~\ref{fig:newislice}).

\begin{figure}
\begin{center}
\includegraphics[width=8cm]{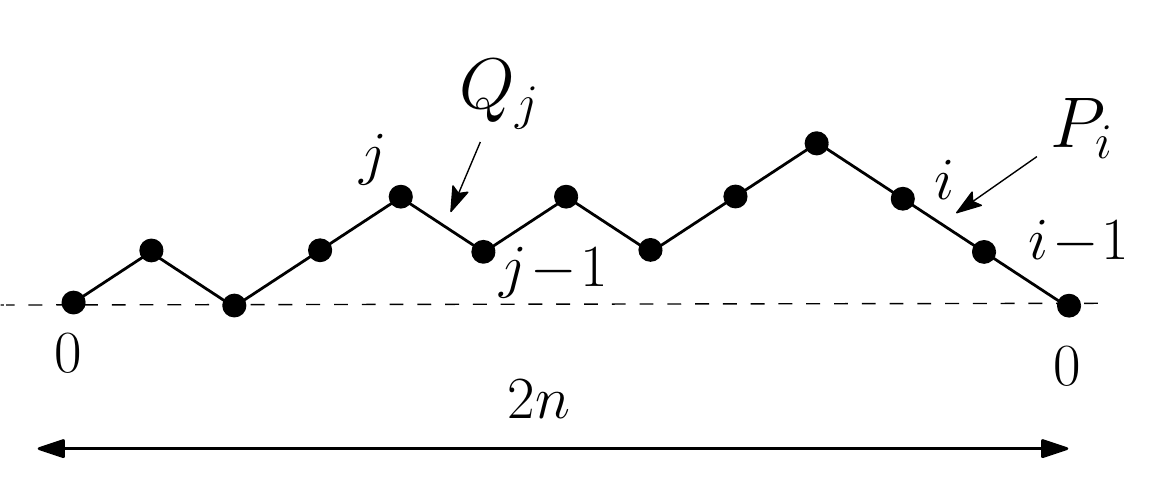}
\end{center}
\caption{An example of directed path of length $2n=12$, made of elementary steps with height difference $\pm 1$, 
starting and ending at height $0$ and remaining (weakly) above height $0$. 
To compute $J_n$, we must sum over such paths
with a weight $Q_i$ (resp.\ $P_i$) assigned to each descending step $i\to i-1$ following an ascent (resp.\ following a descent).}
\label{fig:pathPiQi}
\end{figure}

Returning to the slice decomposition, it now implies, for any positive $n$, that
\begin{equation*}
J_n=\hat{Z}_{0,0}^{+}(2n;\{P_i\}_{i\geq 1},\{Q_i\}_{i\geq 1})\ ,
\end{equation*}
where $\hat{Z}_{0,0}^{+}(2n;\{P_i\}_{i\geq 1},\{Q_i\}_{i\geq 1}$ denotes the generating function of paths of length $2n$,
made of elementary steps of height difference $\pm 1$, starting and ending at height $0$ and remaining above height $0$, with each descending step from height 
$i$ to height $i-1$ weighted by $P_i$ if it follows a descending step $i+1\to i$ and by $Q_i$ if it follows an ascending step $i-1\to i$
(see Fig.~\ref{fig:pathPiQi}).
Setting $J_0=1$ and using the shorthand notation $J(z)=J(z;t_\bullet,t_\circ)$, this is summarized into the new continued fraction expansion 
\begin{equation}
\begin{split}
J(z)&\equiv \sum_{n\geq 0} J_n z^n\\
& =\frac{1}{\displaystyle{1\!-\!z (Q_1\!-\!P_1)\!-z \frac{P_1}{\displaystyle{1\!-\!z (Q_2\!-\!P_2)\!-\!z \frac{P_2}{\displaystyle{1\!-\!z (Q_3\!-\!P_3)\!-\!z \frac{P_3}{\displaystyle{ 1\!-\! \cdots}}}}}}}}\\
\end{split}
\label{eq:contfracnew}
\end{equation}
which we may write as
\begin{equation}
J(z)=\frac{1}{\displaystyle{1-z Y_1-z \frac{Y_2}{\displaystyle{1-z Y_3-z \frac{Y_4}{\displaystyle{1-z Y_5-z \frac{Y_6}{\displaystyle{ 1- \cdots}}}}}}}}
\label{eq:contfracnewbis}
\end{equation}
upon defining
\begin{equation}
Y_{2i-1}\equiv Q_i-P_i\ , \qquad Y_{2i}= P_i
\label{eq:defYi}
\end{equation}
for $i\geq 1$.
\begin{figure}
\begin{center}
\includegraphics[width=8cm]{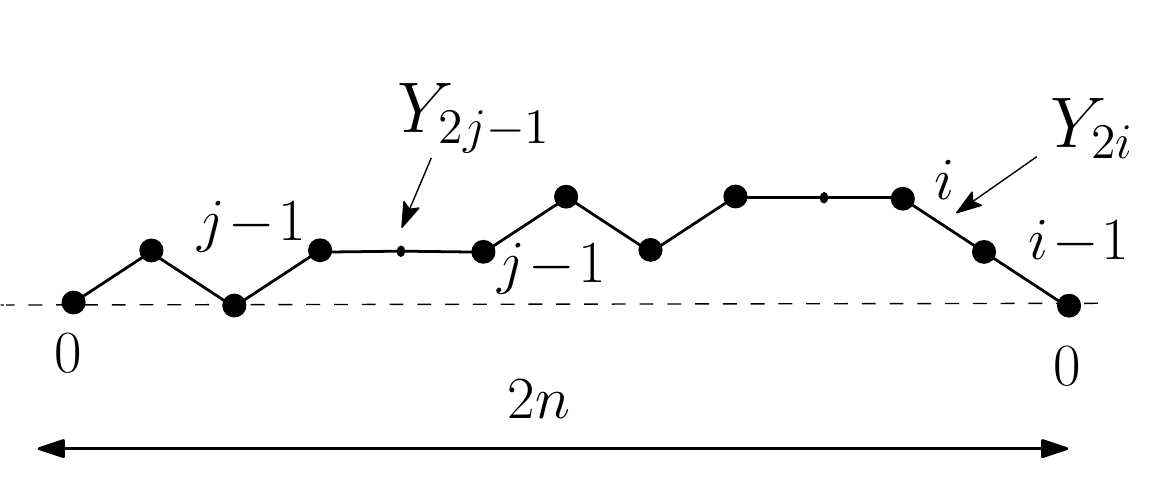}
\end{center}
\caption{An example of directed path of length $2n=12$, made of elementary steps with height difference $\pm 1$
and elongated steps (of horizontal length $2$) with height difference $0$,
starting and ending at height $0$ and remaining (weakly) above height $0$. 
To compute $J_n$, we must sum over such paths
with a weight $Y_{2i-1}$ (resp.\ $Y_{2i}$) assigned to each elongated step $i-1\to i-1$ (resp.\ each descending step $i\to i-1$).}
\label{fig:pathYi}
\end{figure}
To understand \eqref{eq:contfracnew}, or equivalently \eqref{eq:contfracnewbis}, we note that, expanding the right hand side of this latter equation, the term of order $z^n$ builds the generating function 
$\check{Z}_{0,0}^+(2n;\{Y_i\}_{i\geq 1})$ of paths of length $2n$ starting and ending at height $0$ and remaining above height $0$,
made of elementary (i.e.\ of horizontal length $1$) steps of height difference $\pm 1$ together with ``elongated steps" of horizontal
length $2$ and height difference $0$. Each elementary descending step from height $i$ to height $i-1$ ($i\geq 1$) receives 
a weight $Y_{2i}$ while each elongated step at height $i-1$ ($i\geq 1$) receives the weight $Y_{2i-1}$ (see Fig.~\ref{fig:pathYi}).
Deforming each elongated step at height $i-1$ into a sequence of elementary steps $i-1\to i\to i-1$, we recover paths 
made only of elementary steps of height difference $\pm 1$, and (after regrouping all paths with the same deformation) receiving a weight $Y_{2i-1}+Y_{2i}=Q_i$ for
each sequence $i-1\to i \to i-1$ or equivalently for each descending step $i\to i-1$ following an ascent, and a weight $Y_{2i}=P_i$ for those elementary steps $i\to i-1$ which are
not part of a sequence $i-1\to i \to i-1$, i.e.\ follow a descent. In other words, 
$\check{Z}_{0,0}^+(2n;\{Y_i\}_{i\geq 1})=\hat{Z}_{0,0}^{+}(2n;\{P_i\}_{i\geq 1},\{Q_i\}_{i\geq 1})$, which explains the identity \eqref{eq:contfracnew}.

To conclude this section, let us rewrite our fundamental equality \eqref{eq:fundequal} in the more compact form
\begin{equation*}
J(z;t_\bullet,t_\circ)=F(z;t_\bullet,t_\circ)\ .
\end{equation*}

\begin{figure}
\begin{center}
\includegraphics[width=3cm]{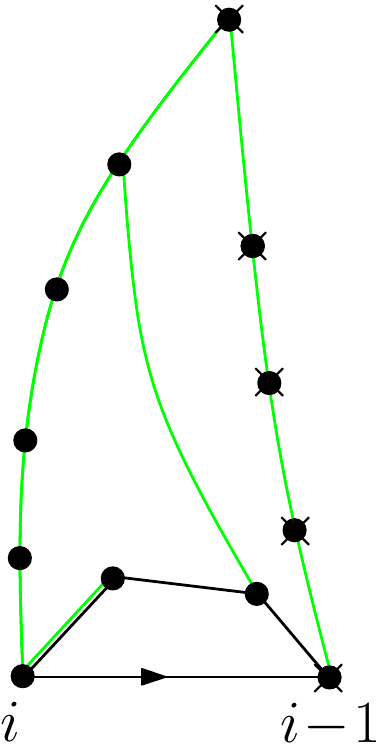}
\end{center}
\caption{A schematic picture of the slice decomposition of a (non-trivial) $i$-slice into two slices, leading
to the recursion relations \eqref{eq:BiWirecur} and \eqref{eq:PiQirecur}.}
\label{fig:recurislice}
\end{figure}

\section{Getting the slice generating functions by solving recursion relations}
\label{sec:recureq}
The slice generating functions $B_i$, $W_i$, $P_i$ and $Q_i$ satisfy systems of non linear recursive equations
which may be derived by performing a slice decomposition of the slices themselves.
Indeed, when the $i$-slice, of left boundary length $\ell$, is not reduced to a single edge, we may look at the sequence of vertices
encountered clockwise around the face lying on the left of the root-edge of the slice and draw the leftmost geodesic 
paths from these vertices to the apex. Using the labeling $d_s(v)+i-\ell$ , 
the sequence of encountered labels, starting from its root-vertex,
is either $i\to i+1\to i\to i-1$, $i\to i-1\to i\to i-1$ or $i\to i-1\to i-2\to i-1$ (if $i\geq 2$)
and, upon cutting along the leftmost geodesic paths, a new slice arises for each descending step of this sequence
(see Fig.~\ref{fig:recurislice}).  

For the first weighting, this decomposition, applied to $i$-slices enumerated by $B_i$ and $W_i$, leads to 
the system
\begin{equation}
 \begin{split}
&B_i=t_\bullet+B_i(W_{i-1}+B_i+W_{i+1})\\
&W_i=t_\circ+W_i(B_{i-1}+W_i+B_{i+1})\\
\end{split}
\label{eq:BiWirecur}
\end{equation}
for $i\geq 1$, with $B_0=W_0=0$. 
For the second weighting, this decomposition, applied to $i$-slices enumerated by $P_i$ and $Q_i$, leads similarly to 
the system
\begin{equation}
 \begin{split}
&P_i=t_\bullet+P_i(P_{i-1}+Q_i+Q_{i+1})\\
&Q_i=t_\circ+Q_i(P_{i-1}+Q_i)+P_iQ_{i+1}\\
\end{split}
\label{eq:PiQirecur}
\end{equation}
for $i\geq 1$, with $P_0=Q_0=0$. 

The solution of \eqref{eq:BiWirecur} was derived in \cite{FG14}. Parametrizing $t_\bullet$ and $t_\circ$ by $x$ and $\gamma$ via
\begin{equation}
\begin{split}
&t_\bullet=\frac{x(\gamma-x)^3(1-\gamma x^3)}{(x+x^3+\gamma-6\, x^2 \gamma+x^4 \gamma+x\, \gamma^2+x^3 \gamma^2)^2}\\
&t_\circ=\frac{x(\gamma-x^3)(1-\gamma x)^3}{(x+x^3+\gamma-6\, x^2 \gamma+x^4 \gamma+x\, \gamma^2+x^3 \gamma^2)^2}\ ,\\
\end{split} 
\label{eq:ttparam}
\end{equation}
with $|x|\leq 1$~\footnote{The parametrization in invariant under $(x,\gamma)\to(1/x,1/\gamma)$ so we may always choose $|x|\leq 1$.}, 
it was shown that 
\begin{equation}
\begin{split}
& B_{2i}=B \frac{(1-x^{2i})(1-\gamma x^{2i+3})}{(1-\gamma x^{2i+1})(1-x^{2i+2}) }\qquad   W_{2i+1}=W\frac{(1-\gamma x^{2i+1})(1-x^{2i+4})}{(1-x^{2i+2})(1-\gamma x^{2i+3})}\\
& B_{2i+1}=B \frac{(1-x^{2i+1}/\gamma)(1-x^{2i+4})}{(1-x^{2i+2})(1-x^{2i+3}/\gamma)}\qquad   W_{2i}=W\frac{(1-x^{2i})(1-x^{2i+3}/\gamma)}{(1-x^{2i+1}/\gamma)(1-x^{2i+2})}\\
\end{split}
\label{eq:BiWiexpl}
\end{equation}
for $i\geq 0$, where
\begin{equation*}
\begin{split}
&B=\frac{x(\gamma-x)^2}{x+x^3+\gamma-6\, x^2 \gamma+x^4 \gamma+x\, \gamma^2+x^3 \gamma^2}\\
&W=\frac{x(1-\gamma x)^2}{x+x^3+\gamma-6\, x^2 \gamma+x^4 \gamma+x\, \gamma^2+x^3 \gamma^2}\ .\\
\end{split} 
\end{equation*}

As for the solution of \eqref{eq:PiQirecur}, it was guessed in \cite{AmBudd}. Parametrizing now $t_\bullet$ and $t_\circ$ by $y$ and $\alpha$ via
\begin{equation}
\begin{split}
&t_\bullet= \frac{y (1 - \alpha y)^3 (1 - \alpha y^3)}{ (1+y+\alpha y-6\, \alpha y^2+\alpha y^3 + \alpha^2 y^3 + \alpha^2 y^4)^2}\\
&t_\circ = \frac{\alpha y (1 - y)^3 (1 - \alpha^2 y^3)}{(1+y+\alpha y-6\, \alpha y^2+\alpha y^3 + \alpha^2 y^3 + \alpha^2 y^4)^2  }\ , \\
  \end{split}
  \label{eq:ttparambis}
\end{equation}
with $|y|\leq 1$, it was found that 
\begin{equation}
P_i = P \frac{(1-y^i)(1-\alpha y^{i+3})}{(1-y^{i+1})(1-\alpha y^{i+2})}
\qquad
Q_i = Q \frac{(1-y^i)(1-\alpha^2 y^{i+3})}{(1-\alpha y^{i+1})(1-\alpha y^{i+2})}\ ,
\label{eq:PiQiexpl}
\end{equation}
for $i\geq 0$, where
\begin{equation}
\begin{split}
&P=\frac{y(1-\alpha y)^2}{1+y+\alpha y-6\, \alpha y^2+\alpha y^3 + \alpha^2 y^3 + \alpha^2 y^4}\\
&Q=\frac{\alpha y(1-y)^2}{1+y+\alpha y-6\, \alpha y^2+\alpha y^3 + \alpha^2 y^3 + \alpha^2 y^4}\ .\\
\end{split} 
\label{eq:PQparam}
\end{equation}
Note that the two parametrizations \eqref{eq:ttparam} and \eqref{eq:ttparambis} are actually equivalent providing we relate $y$ and $\alpha$ 
to $x$ and $\gamma$ via
\begin{equation*}
\alpha=\frac{1}{\gamma^2}\ , \qquad y=\gamma x\ .
\end{equation*} 
With this correspondence, we immediately deduce that
\begin{equation*}
P=B\ , \qquad Q=W\ .
\end{equation*}
This should not come as a surprise since, from  \eqref{eq:BiWiexpl} and \eqref{eq:PiQiexpl}, $B$, $W$, $P$ and $Q$ are the $i\to \infty$ limits of $B_i$, $W_i$, $P_i$ and $Q_i$, enumerating slices \emph{with no bound on their boundary lengths}. From \eqref{eq:BiWirecur} and \eqref{eq:PiQirecur},
both pairs $(B,W)$, and $(P,Q)$ are determined by the \emph{same closed system}, namely:
\begin{equation}
\begin{split}
& B=t_\bullet+B(B+2W)\ , \qquad W= t_\circ+W(W+2B)\ , \\
& P=t_\bullet+P(P+2Q)\ , \qquad Q= t_\circ+Q(Q+2P) \ .\\
\end{split}
\label{eq:BWPQ}
\end{equation}  
Let us end this section by rewriting the results for $P_i$ and $Q_i$ in terms of $Y_i$, as defined in \eqref{eq:defYi}. First, Eq.~\eqref{eq:PiQirecur}
may be rewritten as
\begin{equation}
 \begin{split}
&Y_{2i}=t_\bullet+Y_{2i}(Y_{2i-2}+Y_{2i-1}+Y_{2i}+Y_{2i+1}+Y_{2i+2})\\
&Y_{2i-1}=(t_\circ-t_\bullet)+Y_{2i-1}(Y_{2i-2}+Y_{2i-1}+Y_{2i})\\
\end{split}
\label{eq:Yirecur}
\end{equation}
for $i\geq 1$, with $Y_0=0$. From \eqref{eq:PiQiexpl}, we immediately deduce the solution
\begin{equation}
Y_{2i} = P \frac{(1-y^i)(1-\alpha y^{i+3})}{(1-y^{i+1})(1-\alpha y^{i+2})}
\qquad
Y_{2i+1}=Y \frac{(1-y^{i+1})(1-\alpha y^{i+3})}{(1-y^{i+2})(1-\alpha y^{i+2})}
\ ,
\label{eq:Yiexpl}
\end{equation}
for $i\geq 0$, with
\begin{equation*}
Y=Q-P\ .
\end{equation*}
The aim of this paper is to go beyond the guessing approach of \cite{AmBudd} and to provide a constructive way to obtain this latter formula \eqref{eq:Yiexpl}, and consequently 
\eqref{eq:PiQiexpl}, upon using general results for continued fractions of the type \eqref{eq:contfracnewbis}.
This is indeed the constructive approach used in \cite{FG14} to  obtain the expressions \eqref{eq:BiWiexpl} from general results 
for continued fractions of the type \eqref{eq:contfrac}.

\section{Getting the slice generating functions by extracting continued fraction coefficients: generalities}
\label{sec:extraction}
\subsection{The Stieltjes type}
\label{sec:stieltjes}
Eq.~\eqref{eq:contfrac} is a continued fraction of the so-called Stieltjes type.
Its coefficients $B_{2i}$ and $W_{2i-1}$ for $i\geq 1$ are known to be related to the coefficients $F_n$
via the relations
\begin{equation}
B_{2i}=\frac{h_i^{(0)}}{h_{i-1}^{(0)}}\Big{/}\frac{h_{i-1}^{(1)}}{h_{i-2}^{(1)}}\qquad\qquad\qquad
W_{2i-1}=\frac{h_{i-1}^{(1)}}{h_{i-2}^{(1)}}\Big{/}\frac{h_{i-1}^{(0)}}{h_{i-2}^{(0)}}
\label{eq:BiWiHankel}
\end{equation}
for $i\geq 1$, in terms of the Hankel determinants
\begin{equation*}
h_i^{(0)}=\det (F_{n+m})_{0\leq n,m \leq i}\qquad\qquad
h_i^{(1)}=\det (F_{n+m+1})_{0\leq n,m \leq i}
\end{equation*}
for $i\geq 0$, with the convention $h_{-1}^{(0)}=h_{-1}^{(1)}=1$. These expressions were used in \cite{FG14}
to obtain the expressions \eqref{eq:BiWiexpl} for $B_{2i}$ and $W_{2i+1}$ ($i\geq 0$).
As for the expressions of $B_{2i+1}$ and $W_{2i}$, it is clear from \eqref{eq:BiWirecur} that $B_i$ and $W_i$ 
play symmetric roles upon exchanging $t_\bullet$ and $t_\circ$. The expressions \eqref{eq:BiWiexpl} for $B_{2i+1}$ and $W_{2i}$ 
are simply deduced upon this transformation, which amounts to a change $\gamma \leftrightarrow 1/\gamma$, 
$B\leftrightarrow W$ in the formulas (see \cite{FG14}).
At this stage, it is important to note that the knowledge of the generating functions $F_n$ is not sufficient to determine
all the $B_i$'s and $W_i$'s as the associated continued fraction involves only one parity of the index $i$ ($B_i$'s with even index $i$ and $W_i$'s 
with odd index $i$) and that we have to rely on a symmetry principle to get the other parity. Otherwise stated,
the derivation of all the $B_i$'s and $W_i$'s requires in principle the knowledge of a second family of generating functions.
In the present case, these generating functions are nothing but those of rooted quadrangulations with a boundary
of length $n$, bicolored in such a way that their root-vertex is white instead of black. 
Of course, by symmetry, those are nothing but the $F_n(t_\circ,t_\bullet)$, $n\geq 1$ and a simple
symmetry argument is sufficient to conclude. 

\newcommand{\bfeqref}[1]{\textup{{\normalfont(\bf{\ref{#1}}}\normalfont)}}
\subsection{The type of Eq.~\bfeqref{eq:contfracnewbis}}
\label{sec:newtype} 
When dealing with a continued fraction of the type of Eq.~\eqref{eq:contfracnewbis}, a first remark 
should be emphasized: the knowledge of $J_n$ is not sufficient to determine the coefficients $Y_i$.
Indeed, expanding in $z$ gives rise to the first equations:
\begin{equation}
\begin{split}
& J_1=(Y_1+Y_2)\\
& J_2=(Y_1+Y_2)^2+Y_2(Y_3+Y_4)\\
& \vdots \\
\end{split}
\label{eq:firsteps}
\end{equation}
and it is easily seen that, at each step, two new $Y_i$'s appear on the right hand side, so that the system is clearly underdetermined.

As shown in \cite{DiFKe10,DiFKe11}, a full determination of the coefficients $Y_i$ requires, in addition to the set of $J_n$ for $n\geq 1$, the
knowledge of $Y_1$ and of a second family of quantities $\tilde{J}_n\equiv \tilde{J}_n(t_\bullet,t_\circ)$, $n\geq 0$, satisfying
\begin{equation}
{\tilde J}(z) \equiv \sum_{n\geq 0} \tilde{J}_n z^n=\frac{1}{\displaystyle{1-z \tilde{Y}_1-z \frac{\tilde{Y}_2}{\displaystyle{1-z \tilde{Y}_3-z \frac{\tilde{Y}_4}{\displaystyle{1-z \tilde{Y}_5-z \frac{\tilde{Y}_6}{\displaystyle{ 1- \cdots}}}}}}}}
\label{eq:tildeJn}
\end{equation}
where we have defined (assuming $Y_i\neq 0$ for all $i\geq 1$)
\begin{equation}
\tilde{Y}_{2i-1}\equiv\frac{1}{Y_{2i-1}}\ , \qquad \tilde{Y}_{2i}\equiv \frac{Y_{2i}}{Y_{2i-1}Y_{2i+1}}
\label{eq:tildeYi}
\end{equation}
for $i\geq 1$.
Expanding in $z$ now gives rise to the first equations:
\begin{equation}
\begin{split}
& \tilde{J}_1=\frac{(Y_2+Y_3)}{Y_1Y_3}\\
& \tilde{J}_2=\frac{(Y_2+Y_3)^2}{(Y_1 Y_3)^2}+\frac{Y_2(Y_4+Y_5)}{Y_1 Y_3^2 Y_5}\\
& \vdots \\
\end{split}
\label{eq:firstepsbis}
\end{equation}
Knowing $Y_1$, the first equation of  \eqref{eq:firsteps} yields $Y_2$, then the first equation of \eqref{eq:firstepsbis} yields $Y_3$, 
the second equation of \eqref{eq:firsteps} yields $Y_4$, and so on. The $Y_i$'s are now fully determined and a compact 
formula may be written as follows: define 
\begin{equation}
j_n\equiv \left\{
\begin{matrix}
1& \hbox{if}\ n=0 \\
& \\
Y_1\, J_{n-1} & \hbox{if}\ n\geq  1 \\
& \\
{\tilde J}_{-n}& \hbox{if}\ n\leq -1\\
\end{matrix}
 \right.
 \label{eq:jn}
\end{equation}
and the Hankel-type determinants
\begin{equation}
H_i^{(0)}=\det (j_{n+m-i-1})_{1\leq n,m \leq i}\qquad\qquad
H_i^{(1)}=\det (j_{n+m-i})_{1\leq n,m \leq i}\ .
\label{eq:Hankeltype}
\end{equation}
Then we have, for $i\geq 1$, the following formulas, reminiscent of \eqref{eq:BiWiHankel},
\begin{equation}
Y_{2i}=\frac{H_{i-1}^{(0)}}{H_{i}^{(0)}}\Big{/}\frac{H_{i}^{(1)}}{H_{i+1}^{(1)}}\qquad\qquad\qquad
Y_{2i-1}=\frac{H_{i}^{(1)}}{H_{i-1}^{(1)}}\Big{/}\frac{H_{i}^{(0)}}{H_{i-1}^{(0)}}
\label{eq:YiHankel}
\end{equation}
with the convention $H_{0}^{(0)}=H_{0}^{(1)}=0$. A proof of these formulas can be found in \cite{DiFKe10,DiFKe11}. We present a slightly simpler proof in the Appendix A below.
 To summarize, when dealing with a continued fraction of the type of Eq.~\eqref{eq:contfracnewbis},
 we may extract the coefficients $Y_i$ if, in addition to $J(z)$, we also know $Y_1$ 
 and $\tilde{J}(z)$. As we shall see in Sect.~\ref{sec:conserved} below, getting a simple expression for $Y_1$ 
 combinatorially may be achieved upon using a so-called conserved quantity. As for $\tilde{J}(z)$,
 we have not been able to obtain it via combinatorial arguments (as opposed to the previous section, we cannot
 rely here on any symmetry principle to get $\tilde{J}_n(t_\bullet,t_\circ)$ from $J_n(t_\bullet,t_\circ)$).
 Without the knowledge of $\tilde{J}(z)$,
 Eq.~\eqref{eq:contfracnewbis}  yields a much weaker system than the recursion equations \eqref{eq:Yirecur}.
 In fact, \emph{any arbitrary choice of $\tilde{J}_n$ will lead, through \eqref{eq:YiHankel}, to a set of
 $Y_i$'s satisfying Eq.~\eqref{eq:contfracnewbis}}, while the actual $Y_i$'s, solution of Eqs.~\eqref{eq:Yirecur},
 correspond to a unique value of the $\tilde{J}_n$'s, to be determined.

As we shall now explain, we may however \emph{conjecture} a simple expression for $\tilde{J}(z)$,
based on an explicit solution of the problem in the case of \emph{finite} continued fractions.
With this conjectured form of $\tilde{J}(z)$, we may then verify that the obtained $Y_i$'s 
precisely match their actual expressions \eqref{eq:Yiexpl} guessed in \cite{AmBudd}.

\subsection{The case of finite continued fraction}
\label{sec:finite}
In this section, let us briefly digress from our combinatorial problem and discuss the case of a \emph{finite} continued fraction. 
More precisely, let 
\begin{equation*}
J(z)\equiv \frac{1}{\displaystyle{1-z Y_1-z \frac{Y_2}{\displaystyle{1-z Y_3-z \frac{Y_4}{\displaystyle{
\frac{\ddots}{\displaystyle{1-z Y_{2\alpha-3}-z \frac{Y_{2\alpha-2}}{\displaystyle{ 1- z\, Y_{2\alpha-1}}}}}}}}}}}
\end{equation*}
where $Y_1,Y_2,\cdots , Y_{2\alpha -1}$ denote \emph{independent indeterminates}. 
We also define 
\begin{equation*}
\begin{split}
\tilde{J}(z)& \equiv \frac{1}{\displaystyle{1-z \tilde{Y}_1-z \frac{\tilde{Y}_2}{\displaystyle{1-z \tilde{Y}_3-z \frac{\tilde{Y}_4}{\displaystyle{
\frac{\ddots}{\displaystyle{1-z \tilde{Y}_{2\alpha-3}-z \frac{\tilde{Y}_{2\alpha-2}}{\displaystyle{ 1- z\, \tilde{Y}_{2\alpha-1}}}}}}}}}}}
\\
\hskip-1.2cm \hbox{with}\ \ \tilde{Y}_{2i-1}=& \frac{1}{Y_{2i-1}}\ \hbox{for}\ 1\leq i\leq\alpha\quad \hbox{and}\ \ \tilde{Y}_{2i}=\frac{Y_{2i}}{Y_{2i-1}Y_{2i+1}}
\ \hbox{for}\ 1\leq i<\alpha\ .\\
\end{split}
\end{equation*}
The rational function $J(z)$ is easily seen to be the 
ratio of a polynomial of degree $\alpha-1$ in $z$ by a polynomial
of degree $\alpha$ in $z$, hence characterized by $(\alpha-1+1)+(\alpha+1)-1-1=2 \alpha-1$ coefficients (the last two $-1$'s
correspond to removing a global factor in both the numerator and the denominator, and ensuring that
$J(0)=1$), depending on the
$2\alpha-1$ indeterminates $Y_1,Y_2,\cdots ,Y_{2\alpha-1}$. In this case, the knowledge of the function $J(z)$ alone therefore entirely determines all the coefficients of the continued fraction. This property may be reconciled with the apparently contrary statement of the previous section
by noting that, in the present case of a finite continued fraction, both $Y_1$ and $\tilde{J}(z)$ (defined via 
\eqref{eq:tildeJn} and \eqref{eq:tildeYi}) can be deduced from $J(z)$. More precisely, we have the following
relations, derived in Appendix A below:
\begin{equation}
\tilde{J}(z)=-\frac{Y_1}{z}J\left(\frac{1}{z}\right)\ , \qquad Y_1=-\frac{1}{\lim\limits_{z\to \infty}z\, J(z)}\ \ .
\label{eq:JJtilde}
\end{equation}
Note that $\tilde{J}(z)$ is also a rational function of $z$ and that the expression for $Y_1$ simply rephrases the
desired property that $\tilde{J}(z)=1+O(z)$. Knowing $J(z)$, $Y_1$ and $\tilde{J}(z)$, we can then deduce the coefficients $j_n$ for all integer $n$ via
their definition \eqref{eq:jn} and get $Y_2,Y_3,\cdots,Y_{2\alpha-1}$ from Eqs.~\eqref{eq:Hankeltype} and \eqref{eq:YiHankel}
(which are also valid in the case of a finite continued fraction -- see Appendix A).
 The relations \eqref{eq:JJtilde} are proved in the Appendix A below.
 
\section{Recovering \eqref{eq:PiQiexpl} from the continued fraction formalism}
\label{sec:recovering}
\subsection{A conjectured expression for \boldmath{$Y_1$} and \boldmath{$\tilde{J}(z)$}}
\label{sec:conjecture}
Returning now to our enumeration problem of $i$-slices with the second weighting, let us \emph{conjecture} 
that, although our continued fraction is now infinite, \emph{the relations \eqref{eq:JJtilde} still
hold} for the particular choice of  $Y_i$ we are interested in, namely the solution of \eqref{eq:Yirecur}.
More precisely, $J(z)$, originally defined as a power series in $z$, is convergent for small enough real $z$ 
(namely $0\leq z < 1/(\sqrt{Q}+\sqrt{P})^2$, see explicit expressions below -- here we assume that $t_\bullet$ and $t_\circ$ are 
small enough positive reals so that $P$ and $Q$ 
are positive reals) but may be analytically continued to large enough real $z$ ($z>1/(\sqrt{Q}-\sqrt{P})^2$). This allows us to
define $J(1/z)$ for small real $z$ (namely $0\leq z<(\sqrt{Q}-\sqrt{P})^2$) and our conjecture is that, in this range, $\tilde{J}(z)$ is obtained via the relation
$\tilde{J}(z)=-(Y_1/z)J(1/z)$ with a value of $Y_1$ adjusted so that $\tilde{J}(0)=1$.
Assuming this property, let us now see if we can then recover the desired expression \eqref{eq:PiQiexpl}, or equivalently \eqref{eq:Yiexpl}. 

Let us start by recalling the expression of $F(z)$, hence $J(z)$. From \cite{FG14}, we know that
\begin{equation}
F_n=\frac{B}{t_\bullet}(1-B-W)\  Z_{0,0}^{+}(2n;B,W)-\frac{B}{t_\bullet} \ Z_{0,0}^{+}(2n+2;B,W) 
\label{eq:Fnexplicit}
\end{equation}
where $Z_{0,0}^{+}(2n;B,W)$ denotes the generating function of paths of length $2n$,
made of elementary steps with height difference $\pm 1$, colored alternatively in black and white, starting and ending at black height 
$0$ and 
remaining (weakly) above height
$0$, with each descending step from a black height to a white height weighted by
$B$ and each descending step from a white height to a black height weighted by
$W$. A derivation of this expression via slices is recalled in Sect.~\ref{sec:conserved} below.

\begin{figure}
\begin{center}
\includegraphics[width=8cm]{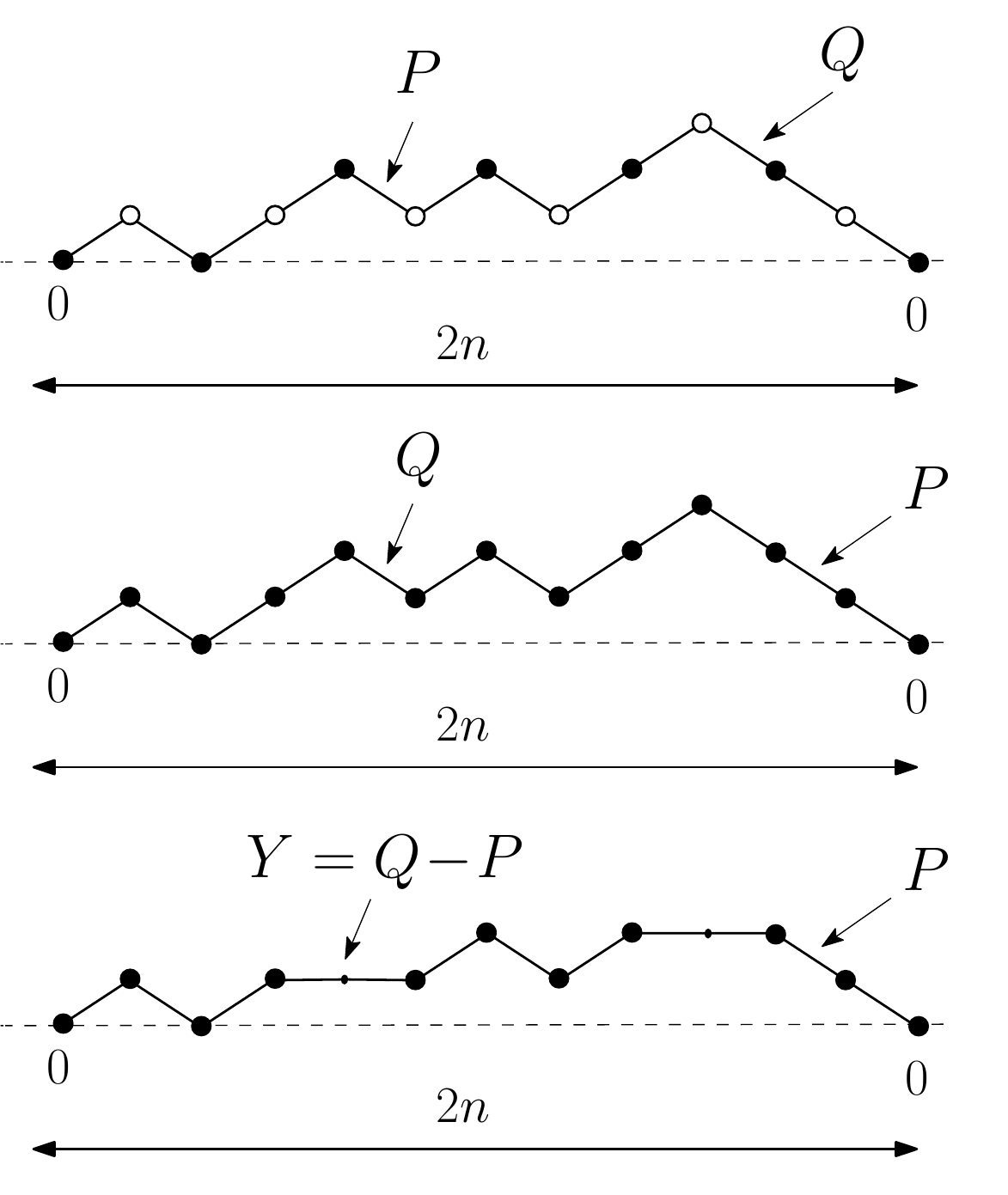}
\end{center}
\caption{A sketch of the three interpretations of the generating function $Z_{0,0}^+(2n;P,Q)$.}
\label{fig:Z2nPQ}
\end{figure}

Equivalently, since $J_n=F_n$, $P=B$ and $Q=W$, we have 
\begin{equation}
\begin{split}
&J_n=A_0\  Z_{0,0}^{+}(2n;P,Q)+A_1\ Z_{0,0}^{+}(2n+2;P,Q) \ ,\\
&A_0=\frac{P}{t_\bullet}(1-P-Q)\ ,\qquad A_1=-\frac{P}{t_\bullet} \ .\\
\end{split}
\label{eq:Jnexpr}
\end{equation}
Let us introduce 
\begin{equation*}
Z(z;P,Q)\equiv\sum_{n\geq 0} Z_{0,0}^{+}(2n;P,Q) z^n\ ,
\end{equation*}
which, by definition, is a solution of
\begin{equation}
Z(z;P,Q)=\frac{1}{\displaystyle{1-z\, \frac{Q}{\displaystyle{1-z\, P\,  Z(z;P,Q)}}}}\ .
\label{eq:eqforZ}
\end{equation}
Note that this (quadratic) equation in $Z$ is equivalent to the equation
\begin{equation*}
Z(z;P,Q)=\frac{1}{\displaystyle{1-z\, (Q-P) -z\,  P\,  Z(z;P,Q)}}
\end{equation*}
so that $Z_{0,0}^{+}(2n;P,Q)$ is also the generating function of paths of length $2n$,
made of elementary steps with height difference $\pm 1$, with each descending step weighted by
$Q$ if it follows an ascending step and by $P$ otherwise.
Alternatively, $Z_{0,0}^{+}(2n;P,Q)$ enumerates paths of length $2n$,
made of elementary steps of horizontal length $1$ and height difference $\pm 1$, and elongated steps
of horizontal length $2$ and height difference $0$, each elongated step receiving the weight $Y=(Q-P)$
and each elementary descending step the weight $P$ (see Fig.~\ref{fig:Z2nPQ}).
As a continued fraction, we thus have
\begin{equation}
Z(z;P,Q)=\frac{1}{\displaystyle{1-z\, Y-z \frac{P}{\displaystyle{1-z\, Y-z \frac{P}{\displaystyle{1-z\, Y-z \frac{P}{\displaystyle{ 1- \cdots}}}}}}}}\ , \qquad Y=Q-P.
\label{eq:Zcondfrac}
\end{equation}
In terms of $Z$, we may write
\begin{equation*}
J(z)= A_0\ Z(z;P,Q)+A_1\  \frac{Z(z;P,Q)-1}{z}
\end{equation*}
and, in components
\begin{equation}
J_n=A_0\ [z^n]Z(z;P,Q)+A_1\ [z^{n+1}]Z(z;P,Q)
\label{eq:JnA0A1}
\end{equation}
for $n\geq 0$.
At this stage, it is important to note that Eq.~\eqref{eq:eqforZ} yields two branches for $Z$ for real $z$, 
namely
\begin{equation*}
\begin{split}
&Z_-(z;P,Q)=\frac{1-Y\, z-\sqrt{(1-Y\, z)^2-4 P\, z}}{2P\, z}\ ,\\
& Z_+(z;P,Q)=\frac{1-Y\, z+\sqrt{(1-Y\, z)^2-4 P\, z}}{2P\, z}\ , \qquad Y=Q-P \\
\end{split}
\end{equation*}
for $|z|\leq 1/(\sqrt{Q}+\sqrt{P})^2$.  
To recover the coefficients $J_n$, we must expand $J(z)$, hence $Z(z;P,Q)$ at small $z$, which requires
to choose 
\begin{equation*}
Z(z;P,Q)=Z_-(z;P,Q).
\end{equation*} 
From $Z_-(z;P,Q)=1+Q z +O(z^2)$ we get $J(z)=1+O(z)$ (since from \eqref{eq:BWPQ}, $A_0+A_1\, Q=1$),
as wanted. 
Using \eqref{eq:JJtilde}, we find that
\begin{equation*}
\begin{split}
&\tilde{J}(z)=-Y_1\left(A_0\, \bar{Z}(z;P,Q)+A_1\, (z\, \bar{Z}(z;P,Q)-1)\right)\\
& \hbox{where}\ \ \bar{Z}(z;P,Q)\equiv \frac{1}{z}Z\left(\frac{1}{z};P,Q\right)\ . \\
\end{split}
\end{equation*}
Again, we have two possible branches for real $z$:
 \begin{equation*}
\begin{split}
&\bar{Z}_-(z;P,Q)=\frac{z-Y-\sqrt{(z-Y)^2-4 P\, z}}{2P\, z}\ ,\\
&\bar{Z}_+(z;P,Q)=\frac{z-Y+\sqrt{(z-Y)^2-4 P\, z}}{2P\, z}\ , \qquad Y=Q-P \\
\end{split}
\end{equation*}
and, to get the $\tilde{J}_n$'s, we must, depending on whether $P>Q$ or $Q>P$, choose the first or second branch respectively
to get rid of the $1/z$ term when $z\to 0$. Both situations  yield actually the same expression for $\tilde{J}_n$.
Assuming for instance $P>Q$, we get 
\begin{equation*}
\tilde{J}_n=-Y_1\, A_0\ [z^n]\bar{Z}_-(z;P,Q)-Y_1\,A_1\ [z^n](z \bar{Z}_-(z;P,Q)-1)\ .
\end{equation*}
The value of $Y_1$ is obtained by ensuring that $\tilde{J}_0=1$. Using $\bar{Z}_-(z;P,Q)=-1/Y+O(z)$, 
we deduce 
\begin{equation}
Y_1\left( \frac{A_0}{Y}+A_1\right)=1\ ,\ \hbox{hence}\ 
Y_1=\frac{Y}{A_0+A_1\, Y}= \frac{(Q-P)(1-P-2Q)}{1-2 Q}\ .
\label{eq:Y1val}
\end{equation}
This value matches that obtained directly via conserved quantities in Sect.~\ref{sec:conserved}.

Using the identity
\begin{equation*}
\bar{Z}_-(z;P,Q)=\frac{1}{P-Q}Z_-\left(\frac{z}{(P-Q)^2};P,Q\right)\ ,
\end{equation*}
we eventually arrive at
\begin{equation}
\begin{split}
\tilde{J}_n&=-\frac{Y_1}{(P-Q)^{2n+1}}\, \left(A_0\ [z^n]Z_-(z;P,Q)+A_1\ (P-Q)^2 [z^n](z Z_-(z;P,Q)-1)\right)\\
&= \frac{Y_1}{(Q-P)^{2n+1}}\, \left(
A_0\  Z_{0,0}^{+}(2n;P,Q)+A_1(Q-P)^2\ Z_{0,0}^{+}(2n-2;P,Q)\right) \\
\end{split}
\label{eq:tildeJnval}
\end{equation}
for $n>0$. For $Q>P$, we must use instead $\bar{Z}_+(z;P,Q)$ but again in this case, $\bar{Z}_+(z;P,Q)=-1/Y+O(z)$
and $\bar{Z}_+(z;P,Q)=\frac{1}{P-Q}Z_-\left(\frac{z}{(P-Q)^2};P,Q\right)$ so that the expressions \eqref{eq:Y1val} and \eqref{eq:tildeJnval} 
remain unchanged.
The first line of \eqref{eq:tildeJnval} may be rewritten as
\begin{equation}
\begin{split}
\tilde{J}_n&=\frac{Y_1}{Y}\, \left(A_0\ [z^n]Z_-\left(\frac{z}{Y^2};P,Q\right)+A_1\ [z^n]\left(z Z_-\left(\frac{z}{Y^2};P,Q\right)-1\right)\right)\\
&= \frac{Y_1}{Y}\,\left(A_0\ [z^n]\tilde{Z}(z;P,Q)+A_1\ [z^{n-1}]\tilde{Z}(z;P,Q)\right)\\
\end{split}
\label{eq:tJnA0A1}
\end{equation}
for $n\geq 1$, with
\begin{equation*}
\begin{split}
&\tilde{Z}(z;P,Q)\equiv Z_-\left(\frac{z}{Y^2};P,Q\right)=\frac{1}{\displaystyle{1-z\, \tilde{Y}-z \frac{\tilde{P}}{\displaystyle{1-z\, \tilde{Y}-z \frac{\tilde{P}}{\displaystyle{1-z\, \tilde{Y}-z \frac{\tilde{P}}{\displaystyle{ 1- \cdots}}}}}}}}\\
&  \hbox{where}\ \tilde{Y}=\frac{1}{Y}\ \hbox{and}\ \tilde{P}=\frac{P}{Y^2}\ .\\
\end{split}
\end{equation*}
Upon defining
\begin{equation*}
k_n\equiv \left\{
\begin{matrix}
1& \hbox{if}\ n=0 \\
& \\
Y\, Z_{n-1} & \hbox{if}\ n\geq  1 \\
& \\
{\tilde Z}_{-n}& \hbox{if}\ n\leq -1\\
\end{matrix}
 \right.
\end{equation*}
where $Z_n\equiv [z^n] Z(z;P,Q)=Z_{0,0}^{+}(2n;P,Q)$ and ${\tilde Z}_{n}\equiv [z^n] \tilde{Z}(z;P,Q)=Z_{0,0}^{+}(2n;P,Q)/Y^{2n}$ ($n\geq 0$),
we may summarize  \eqref{eq:JnA0A1}, \eqref{eq:Y1val} and \eqref{eq:tJnA0A1}  into
\begin{equation*}
j_n=A_0\ \frac{Y_1}{Y}\ k_n+A_1\ \frac{Y_1}{Y}\ k_{n+1} \ \hbox{for all integer}\ n.
\end{equation*}
\subsection{Computation of \boldmath{$H_i^{(0)}$} and \boldmath{$H_i^{(1)}$}}
\label{sec:comput}
The above expression for $j_n$ for all integer $n$ 
opens the way to compute $H_i^{(0)}$ and $H_i^{(1)}$ via \eqref{eq:Hankeltype}.
Indeed, as shown in Appendix B, the coefficients $k_n$ satisfy a set of linear relations of the form
\begin{equation}
\sum_{m=0}^{i-1}k_{n-i+m} (-1)^{m}\frac{x^{(i-1)}_{i-1-m}}{x^{(i-1)}_{i-1}}=0\ , \qquad 2\leq n\leq i
\label{eq:linearA}
\end{equation} 
while, for $n=1$ and $n=i+1$, we have
\begin{equation}
\begin{split}
&\sum_{m=0}^{i-1}k_{1-i+m} (-1)^{m}\frac{x^{(i-1)}_{i-1-m}}{x^{(i-1)}_{i-1}}=\frac{P^{i-1}}{Y^{2(i-1)}}\\
&\sum_{m=0}^{i-1}k_{1+m}(-1)^{m}\frac{x^{(i-1)}_{i-1-m}}{x^{(i-1)}_{i-1}}=(-1)^{i-1}\frac{P^{i-1}}{Y^{i-2}}\ .\\
\end{split}
\label{eq:linearB}
\end{equation}
From these relations, replacing the first column $(j_{n-i})_{1\leq n\leq i}$ of $H_i^{(0)}$ by the linear combination
$\left(\sum\limits_{m=0}^{i-1}j_{n-i+m} (-1)^{m}{x^{(i-1)}_{i-1-m}}/{x^{(i-1)}_{i-1}}\right)_{1\leq n\leq i}$ of this first column with the $i-1$ 
last ones allows us to write
\def\rddots{.^{\displaystyle{\,.}^{\displaystyle{\,.}}}} 
\begin{equation*}
\hskip -1.2cm H_i^{(0)}= \left\vert 
\begin{matrix}
\displaystyle{A_0 \frac{Y_1}{Y} \frac{P^{i-1}}{Y^{2(i-1)}}} & j_{-i+2} & \cdots & \cdots & j_0 \\
& & & &\\
0 & j_{-i+3} & \cdots  & \cdots & j_1\\
& & & &\\
\vdots & \vdots & \rddots & \rddots & \vdots \\
& & & & \\
0 & j_0 & \rddots  & \rddots &  j_{i-2} \\
& & & &\\
\displaystyle{(-1)^{i-1} A_1 \frac{Y_1}{Y} \frac{P^{i-1}}{Y^{i-2}} }& j_1 & \cdots & \cdots & j_{i-1} \\
\end{matrix} \right\vert
= A_0 Y_1 \frac{P^{i-1}}{Y^{2i-1}}H_{i-1}^{(1)}+A_1 Y_1 \frac{P^{i-1}}{Y^{i-1}}H_{i-1}^{(0)}\ .
\end{equation*}
Alternatively, the coefficients $k_n$ satisfy another set of linear relations of the form (see Appendix B)
\begin{equation}
\sum_{m=0}^{i-1}k_{n-m} (-1)^{m}x'^{(i-1)}_{m}=0\ , \qquad 2\leq n\leq i
\label{eq:linearC}
\end{equation} 
while, for $n=1$ and $n=i+1$, we have
\begin{equation}
\begin{split}
&\sum_{m=0}^{i-1}k_{1-m} (-1)^{m}x'^{(i-1)}_{m}=(-1)^{i-1} \frac{P^{i-1}}{Y^{i-2}}\\
&\sum_{m=0}^{i-1}k_{i+1-m}(-1)^{m}x'^{(i-1)}_{m}=Y\, P^{i-1}(Y+P)\ .\\
\end{split}
\label{eq:linearD}
\end{equation}
From these relations, replacing the last column $(j_{n})_{1\leq n\leq i}$ of $H_i^{(1)}$ by the linear combination
$\left(\sum\limits_{m=0}^{i-1}j_{n-m} (-1)^{m}x'^{(i-1)}_{m}\right)_{1\leq n\leq i}$ of this last column with the $i-1$ 
first ones, and using $x'^{(i-1)}_0=1$ (see Appendix B), allows us to write
\begin{equation*}
\hskip -1.2cm H_i^{(1)}\!=\! \left\vert 
\begin{matrix}
j_{-i+2} \!& \!\cdots \!&\! \cdots \!&\! j_0\! & \! \displaystyle{(-1)^{i-1} A_0 \frac{Y_1}{Y} \frac{P^{i-1}}{Y^{i-2}} }\!\\
& & & &\\
 j_{-i+3} \!&\! \cdots \! &\! \cdots \!&\! j_1 \!&\! 0\!\\
& & & &\\
\vdots &\! \vdots \!&\! \rddots \!&\! \rddots \!& \!\vdots \!\\
& & & & \\
 j_0 \!& \!\rddots \! &\! \rddots \!& \! j_{i-2} \!& \!0\! \\
& & & &\\
 j_1 \!& \!\cdots \!&\! \cdots \!&\! j_{i-1}\!&\! \displaystyle{ A_1 \frac{Y_1}{Y}  Y P^{i-1}(Y\!+\!P)}\!\\
\end{matrix} \right\vert
\!=\! A_0 Y_1 \frac{P^{i-1}}{Y^{i-1}}H_{i-1}^{(1)}+A_1 Y_1 P^{i-1}(Y+P)H_{i-1}^{(0)}\ .
\end{equation*}

\begin{figure}
\begin{center}
\includegraphics[width=6cm]{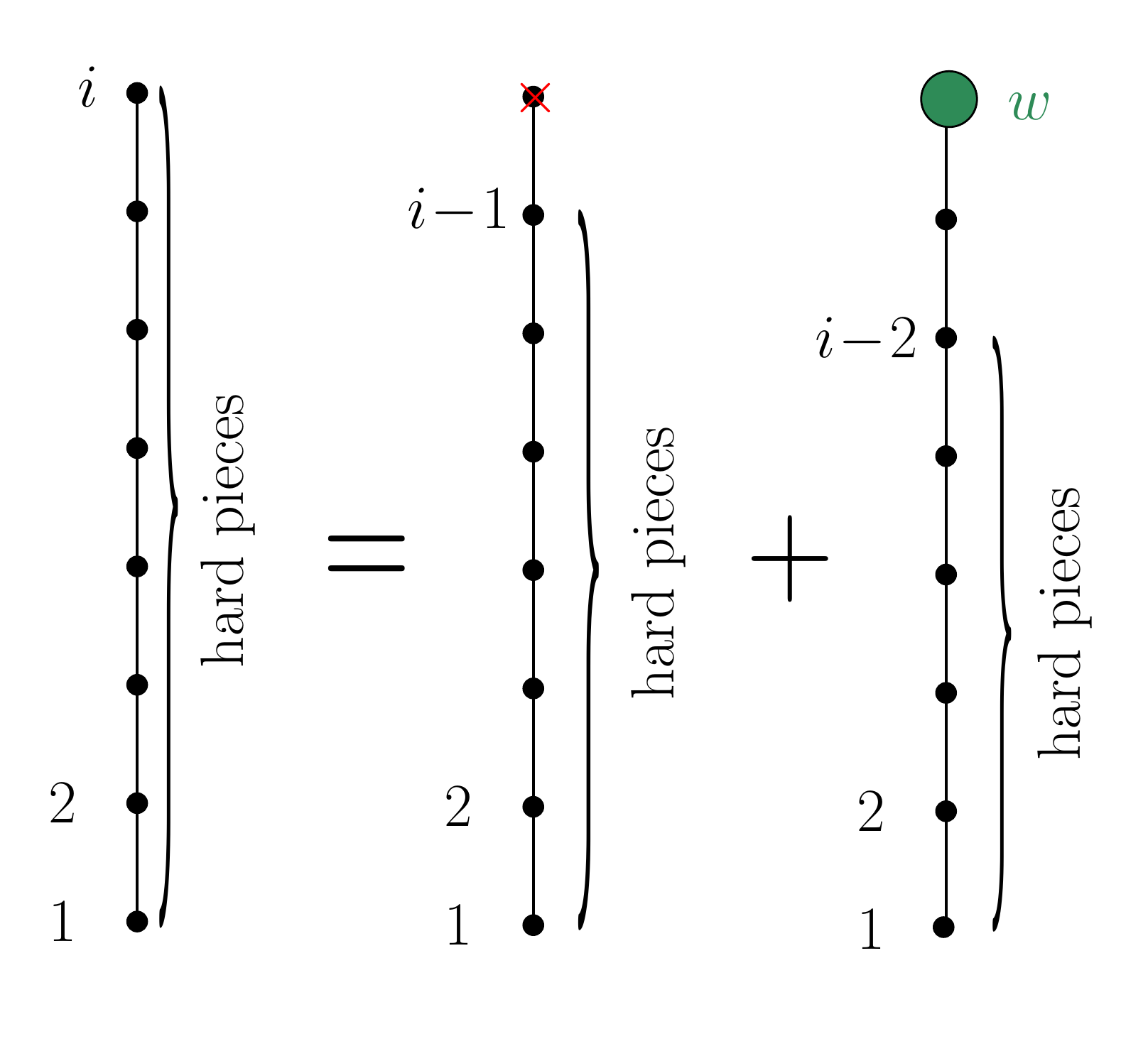}
\end{center}
\caption{A schematic picture of Eq.~\eqref{eq:recurhardpieces}, identifying $L_i^{(0)}$ as the generating function of
hard pieces on a linear graph with $i-1$ vertices.}
\label{fig:hardparticles}
\end{figure}

To summarize,  $H_i^{(0)}$ and $ H_i^{(1)}$ are fully determined by the system 
\begin{equation*}
\begin{split}
& H_i^{(0)}=  A_0 Y_1 \frac{P^{i-1}}{Y^{2i-1}}H_{i-1}^{(1)}+A_1 Y_1 \frac{P^{i-1}}{Y^{i-1}}H_{i-1}^{(0)}\\
&H_i^{(1)}=A_0 Y_1 \frac{P^{i-1}}{Y^{i-1}}H_{i-1}^{(1)}+A_1 Y_1 P^{i-1}(Y+P)H_{i-1}^{(0)}\\
\end{split}
\end{equation*}
for $i\geq 1$ with $H_0^{(0)}=H_0^{(1)}=1$.
Upon setting
\begin{equation}
L^{(0)}_i\equiv \left(\frac{Y}{P}\right)^{\frac{i(i-1)}{2}}H_i^{(0)}\ , \qquad L^{(1)}_i\equiv \frac{1}{Y^{i-1}} \left(\frac{Y}{P}\right)^{\frac{i(i-1)}{2}}H_i^{(1)}\ ,
\label{eq:Lidef}
\end{equation}
these equations read
\begin{equation}
\begin{split}
& L_i^{(0)}=  A_0 \frac{Y_1}{Y^2}L_{i-1}^{(1)}+A_1 Y_1\, L_{i-1}^{(0)}\\
&L_i^{(1)}=A_0 \frac{Y_1}{Y} L_{i-1}^{(1)}+A_1 Y_1(Y+P)\, L_{i-1}^{(0)}\ .\\
\end{split}
\label{eq:systemL}
\end{equation}
Using the first line to express $L_i^{(1)}$ in terms of $L_i^{(0)}$ and re-injecting the result in the second line
yields an equation for $L_i^{(0)}$ only, namely
\begin{equation*}
 L_{i+1}^{(0)}=Y_1\left(\frac{A_0}{Y}+A_1\right)  L_i^{(0)} +A_0 A_1\frac{Y_1^2}{Y^2}P\,  L_{i-1}^{(0)}
\end{equation*}
for $i\geq 1$ with $L_{0}^{(0)}=1$ and $L_{1}^{(0)}=H_{1}^{(0)}=Y_1\left(\frac{A_0}{Y}+A_1\right)$.
Using \eqref{eq:Y1val} and setting
\begin{equation*}
w\equiv A_0 A_1 \frac{Y_1^2}{Y^2}\, P\ ,
\end{equation*}
we recover the well-known equation 
\begin{equation}
 L_{i+1}^{(0)}= L_i^{(0)} +w\, L_{i-1}^{(0)}
 \label{eq:recurhardpieces}
\end{equation}
for $i\geq 1$ with $L_{0}^{(0)}=L_{1}^{(0)}=1$, which allows us to interpret $L_i^{(0)}$ as the generating function of
hard pieces on a linear graph with $i-1$ vertices (see Fig.~\ref{fig:hardparticles}), with a weight $w$ per piece.
The solution of this equation is known to be (see for instance \cite{BG12} Eq.~(6.11))
\begin{equation}
L_i^{(0)}=\frac{1}{(1+y)^{i}}\ \frac{1-y^{i+1}}{1-y}\ \hbox{where}\ w=-\frac{1}{y+y^{-1}+2}\ .
\label{eq:Li0val}
\end{equation}
If we instead eliminate $L_i^{(0)}$ from the system \eqref{eq:systemL}, we obtain for $L_i^{(1)}$ the very same equation 
\begin{equation*}
 L_{i+1}^{(1)}= L_i^{(1)} +w\, L_{i-1}^{(1)}
\end{equation*}
for $i\geq 1$, now with the initial conditions $L_0^{(1)}=Y$ and $L_1^{(1)}=H_1^{(1)}=j_1=Y_1$ (this value
can also be read from \eqref{eq:systemL} as it yields $L_1^{(1)}=Y_1(A_0+A_1(Y+P))=Y_1(A_0+A_1 Q)=Y_1$).
We immediately deduce 
\begin{equation*}
L_i^{(1)} =Y\, L_i^{(0)} +(Y_1-Y) L_{i-1}^{(0)} 
\end{equation*}
with the convention $L_{-1}^{(0)}=0$. Using \eqref{eq:Li0val},  we obtain
\begin{equation}
\begin{split}
L_i^{(1)}& =\frac{1}{(1+y)^{i}}\ \frac{1}{1-y}\left( Y\, (1-y^{i+1}) +(Y_1-Y) (1-y^{i})(1+y)\right)\\
&= \frac{1}{(1+y)^{i}} Y_1 (1+ d\, y)\left(\frac{1-\alpha\, y^{i+2}}{1-y}\right)\\
& \hbox{ where}\ \ \  d\equiv \frac{Y_1-Y}{Y_1}\ \ \ \hbox{and} \ \ \  \alpha \equiv \frac{1}{y^2}\, \frac{d+y}{1+d y}\ .\\
\end{split}
\label{eq:Li1val}
\end{equation}

\subsection{Comparison with formulas \bfeqref{eq:PiQiexpl}}
Combining \eqref{eq:Lidef} and the explicit values \eqref{eq:Li0val} and \eqref{eq:Li1val}, we obtain from \eqref{eq:YiHankel}
the desired expressions \eqref{eq:Yiexpl}. It simply remains to show that our definitions for $y$ and $\alpha$ of Sect.~\ref{sec:comput},
given by 
\eqref{eq:Li0val} and \eqref{eq:Li1val} just above, match their definitions of Sect.~\ref{sec:recureq} given by \eqref{eq:ttparambis}, or equivalently  \eqref{eq:PQparam} and \eqref{eq:BWPQ}. If so,
\eqref{eq:Yiexpl} is equivalent to \eqref{eq:PiQiexpl} and we are done.

Using for $y$ and $\alpha$ their definitions of Sect.~\ref{sec:recureq} (through \eqref{eq:PQparam} and \eqref{eq:BWPQ}), we obtain for $Y=Q-P$ and $Y_1$ (whose value in terms of $P$ and $Q$ is given by \eqref{eq:Y1val}) the parametrizations 
\begin{equation*}
\begin{split}
&Y=\frac{(\alpha-1)y(1-\alpha y^2)}{1+y+\alpha y-6\, \alpha y^2+\alpha y^3 + \alpha^2 y^3 + \alpha^2 y^4}\\
&Y_1=\frac{(\alpha-1)y(1-\alpha y^3)}{(1+y)(1+y+\alpha y-6\, \alpha y^2+\alpha y^3 + \alpha^2 y^3 + \alpha^2 y^4)}\ .\\
\end{split} 
\end{equation*}
so that 
\begin{equation*}
d=-\frac{y(1-\alpha y)}{1-\alpha y^3}
\ \ \ \hbox{and} \ \ \ 
\frac{1}{y^2}\, \frac{d+y}{1+d y}=\alpha
\end{equation*}
as wanted to match the definition \eqref{eq:Li1val} of $\alpha$ of Sect.~\ref{sec:comput}. 

As for $y$, we use the expressions \eqref{eq:Jnexpr} of $A_0$ and $A_1$ to get the 
parametrizations
\begin{equation*}
\begin{split}
&A_0=\frac{(1-\alpha y^2)^2}{(1-\alpha y)(1-\alpha y^3)}\\
&A_1=-\frac{1+y+\alpha y-6\, \alpha y^2+\alpha y^3 + \alpha^2 y^3 + \alpha^2 y^4}{(1-\alpha y)(1-\alpha y^3)}\ .\\
\end{split} 
\end{equation*}
so that 
\begin{equation*}
w=A_0 A_1 \frac{Y_1^2}{Y^2} P =-\frac{1}{y+y^{-1}+2}
\end{equation*}
as wanted to match the definition \eqref{eq:Li0val} of $y$ of Sect.~\ref{sec:comput}. 

A more constructive approach consists in starting instead from the definitions of $y$ and $\alpha$ of Sect.~\ref{sec:comput}
(through \eqref{eq:Li0val} and \eqref{eq:Li1val}) and recovering the parametrization \eqref{eq:PQparam} of Sect.~\ref{sec:recureq}. From \eqref{eq:Jnexpr} (and $t_\bullet=P(1-P-2Q)$),
$A_0$ can be expressed in terms of $P$ and $Q$, as well as $Y$ ($=Q-P$) and $Y_1$ via \eqref{eq:Y1val}.
This leads to
\begin{equation*}
w\equiv A_0 A_1 \frac{Y_1^2}{Y^2} P =-\frac{P(1-Q-P)}{(1-2Q)^2}=-\frac{1}{y+y^{-1}+2}
\end{equation*}
hence we deduce the (so-called characteristic) equation
\begin{equation*}
(1-2 Q)^2-(2+y+y^{-1})P(1-P-Q)=0\ .
\end{equation*} 
This in turn leads to
\begin{equation*}
d\equiv \frac{Y_1-Y}{Y_1}=-\frac{P}{1-P-2Q}\ , \qquad \alpha\equiv \frac{1}{y^2}\frac{d+y}{1+d y}
=-\frac{P-y(1-P-2Q)}{y^2(1-P-2Q -y P)}\ .
\end{equation*}
Using this latter equation to express $Q$ in terms of $P$, $y$ and $\alpha$,
namely
\begin{equation*}
Q=-\frac{P(1+y)(1-\alpha y^2)}{2 y (1-\alpha y)}+\frac{1}{2}\ ,
\end{equation*}
and plugging this
value in the characteristic equation above, we find that $P$ is determined by
\begin{equation*}
-y(1-\alpha y)^2+P(1+y+\alpha y-6\, \alpha y^2+\alpha y^3 + \alpha^2 y^3 + \alpha^2 y^4)=0\ ,
\end{equation*}
from which \eqref{eq:PQparam} follows.

\section{Conserved quantities}
\label{sec:conserved}

The explicit formulas \eqref{eq:BiWiexpl} (resp.\ \eqref{eq:PiQiexpl} or equivalently \eqref{eq:Yiexpl}) are typical expressions
for the solutions of \emph{discrete integrable systems}. A deeper characterization of the integrability of the system  \eqref{eq:BiWirecur} (resp.\ \eqref{eq:PiQirecur} or \eqref{eq:Yirecur}) is the existence of a number of \emph{discrete conserved quantities}, i.e.\ quantities 
whose expression depends explicitly on some positive integer (called $d$ below) but whose value turns out to be independent of this integer. 
In the case of bicolored quadrangulations, it has already been recognized that these conserved quantities may be obtained  
by looking for a direct combinatorial derivation of $F_n$ in the slice formalism. Let us first recall this construction and then see how it
extends to the case of our second weighting governed by local maxima.

\subsection{Conserved quantities for the first weighting}
The slice decomposition described in~\cite{BG12,FG14} applies more generally to
\emph{pointed rooted} quadrangulations with boundaries, that is, 
quadrangulations with a boundary, having a root-edge on the boundary oriented
with the boundary-face on its right, and having a pointed vertex (which 
might not be incident to the boundary-face); the \emph{label} $d(v)$ of each 
vertex $v$ is now the distance from $v$ to the pointed vertex $v_0$.  
For $Q$ such a map, the \emph{canonical bicoloration} of $Q$ is the vertex bicoloration in black and white 
where the root-vertex (origin of the root) is black and any two adjacent vertices have different colors. As before, a \emph{local maximum} 
(or ``local max" for short) for the distance is a vertex $v$ such that 
$d(v)=d(v')+1$ for every neighbor $v'$ of $v$.  
For $n\geq 1$ and $d\geq 0$, 
let $\cBnd$ be the family of admissible pointed rooted quadrangulations with a boundary
of length $2n$, where the root-vertex is at distance \emph{at most} $d$  
from $v_0$, and is \emph{one} (possibly not unique) \emph{of the boundary-vertices that reach 
the smallest distance from $v_0$}.   
Let $F_n^{(d)}$ be the generating function of $\cBnd$
where each black vertex (resp.\ white vertex) receives weight $\tb$ (resp.\ $\tw$)
except for the pointed vertex that receives weight $1$. 
And let $Z_{d,d}^+(2n;\{B_i\}_{i\geq 1},\{W_i\}_{i\geq 1})$ be the generating
functions of paths of length $2n$ starting and ending at height $d$ and staying
at height at least $d$ all along, made of elementary steps with height difference
$\pm 1$, with each descending step from height $i$ to height $i-1$ weighted
by $B_i$ if $i\equiv d\ \mathrm{mod}\ 2$ and weighted by $W_i$ if $i\equiv d+1\ \mathrm{mod}\ 2$.

\begin{figure}
\begin{center}
\includegraphics[width=10cm]{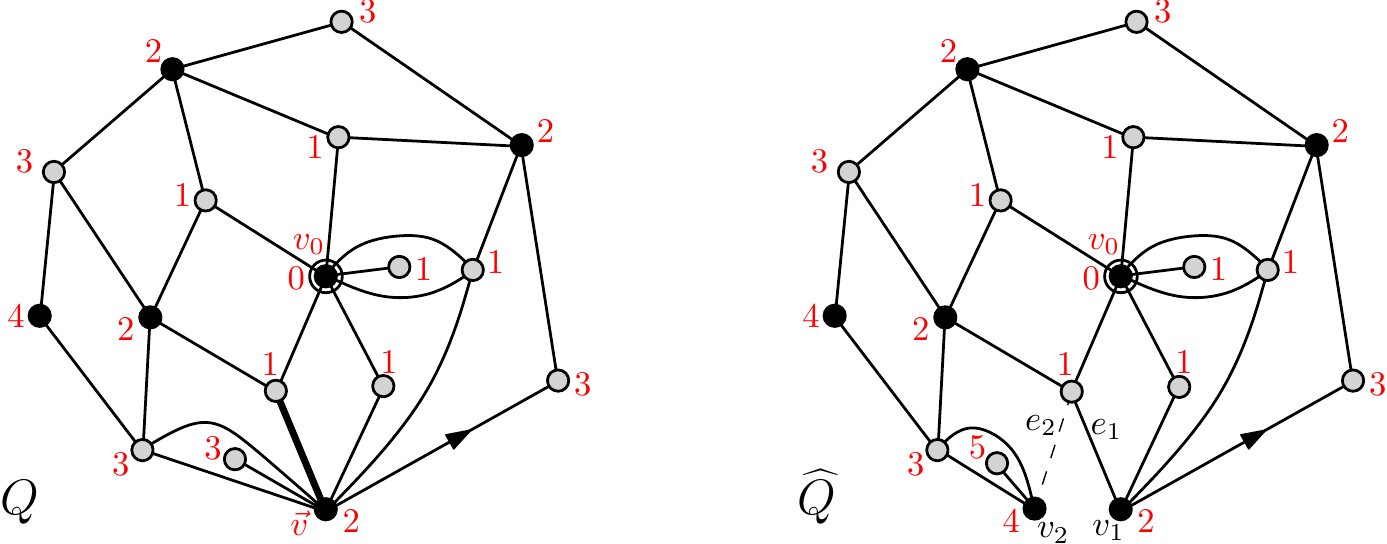}
\end{center}
\caption{Left: a pointed rooted quadrangulation $Q$ with a boundary. Right: the pointed rooted quadrangulation $\hQ$ obtained by cutting along the first
edge of the leftmost geodesic path from the root-vertex to the pointed vertex.}
\label{fig:slitedgequad}
\end{figure}

Note that $\cB_n^{(0)}$ is nothing but the set of rooted quadrangulations with a boundary
of length $2n$ so that $F_n^{(0)}=F_n$. Then, as explained in~\cite{BG12,FG14}, 
the slice decomposition described in Sect.~\ref{sec:contfracFn} for maps in $\cB_n^{(0)}$
applies more generally for maps in $\cB_n^{(d)}$ and yields
$$
F_n^{(d)}=Z_{d,d}^+(2n;\{B_i\}_{i\geq 1},\{W_i\}_{i\geq 1}).
$$
Now let $\cBhnd$ be the subfamily of $\cBnd$ where the pointed vertex is \emph{different from the root-vertex}, and let $\wF_n^{(d)}$ be the generating function for the 
subfamily $\cBhnd$ where the weights are specified as in $F_n^{(d)}$. 
For $Q\in\cBhnd$, with $v_0$ the pointed vertex and $\vec{v}$ the root-vertex,
let $e$ be the first edge of the \emph{leftmost} geodesic path from $\vec{v}$ to $v_0$. 
This edge cannot be a boundary-edge as otherwise, $\vec{v}$ would not reach the  smallest distance from $v_0$
among boundary-vertices.
We may cut along $e$ (starting from $\vec{v}$) so as to duplicate $e$ into two edges
$e_1,e_2$ (with $e_2$ before $e_1$ in ccw order around the new map) and duplicate
$\vec{v}$ into two vertices $v_1,v_2$ (see Fig.\ref{fig:slitedgequad}). Let $\hQ$ be the pointed rooted quadrangulation with a boundary of length $2n+4$ that is  obtained by erasing $e_2$, 
taking $v_1$ as the new root-vertex, and keeping $v_0$ as the pointed vertex. 
Denoting by 
$\mathbf{d}=(d_1,\ldots,d_{2n+4})$ the distances from $v_0$ of the successive boundary-vertices (starting from $v_1$) in ccw order
around $\hQ$, we have the conditions that $d_{i+1}=d_i\pm 1$ for $i\in\{1,\cdots,2n+3\}$, 
$d_1$ equals the distance of $\vec{v}$ from $v_0$ in $Q$ 
so that $d_1 \leq d$, $d_i\geq d_1$ for all $i\in\{1,\cdots,2n+1\}$, $d_{2n+1}>d_1$ (indeed,
by the effect
of cutting along the first edge of the leftmost geodesic path, the distance of $v_2$
from $v_0$ is strictly larger than the distance of $v_1$ from $v_0$), and $d_{2n+4}=d_1-1$. The bipartiteness of $\widehat{Q}$ implies that $d_{2n+1}\equiv d_1\ \mathrm{mod}\ 2$, so that the last entries of $\mathbf{d}$ 
must be $(d_1+2,d_1+1,d_1,d_1-1)$. Hence, if for $k\geq 1$ we denote by 
$Z_{d,d}^{+,k\searrow}(2n;\{B_i\}_{i\geq 1},\{W_i\}_{i\geq 1})$ the generating 
function defined as $Z_{d,d}^+(2n;\{B_i\}_{i\geq 1},\{W_i\}_{i\geq 1})$, but with 
the restriction that the $k$ last steps of the path are descending, then the slice
decomposition applied to $\widehat{Q}$ gives
$$
\hF_n^{(d)}=\frac1{\tb}Z_{d,d}^{+,2\searrow}(2n+2;\{B_i\}_{i\geq 1},\{W_i\}_{i\geq 1})\cdot B_d,
$$
where the factor $\frac1{\tb}$ accounts for the (black) root-vertex being duplicated 
and the factor $B_d$ accounts for the last descent $d_1,d_1-1$.  
Hence for each $n\geq 1$ we have the conserved quantity
\begin{equation}\label{eq:invariantFn}
F_n=Z_{d,d}^+(2n;\{B_i\}_{i\geq 1},\{W_i\}_{i\geq 1})-\frac1{\tb}Z_{d,d}^{+,2\searrow}(2n+2;\{B_i\}_{i\geq 1},\{W_i\}_{i\geq 1})\cdot B_d.
\end{equation}
The first two conserved quantities, $n\in\{1,2\}$, are (with $i=d+1$): for all $i\geq 1$ (with $B_0=0$)
\begin{eqnarray*}
F_1&=&W_i-\frac1{\tb}B_{i+1}W_iB_{i-1},\\
F_2&=&W_i^2+B_{i+1}W_i-\frac1{\tb}(W_i+B_{i+1}+W_{i+2})B_{i+1}W_iB_{i-1}.
\end{eqnarray*}
Shifting in \eqref{eq:invariantFn} all path heights by $-d$ and replacing $B_i$ and $W_i$ by $B_{i+d}$
and $W_{i+d}$ so as to compensate this shift, we get, upon sending $d\to\infty$ the identity
\begin{equation}
F_n=Z_{0,0}^+(2n;B,W)-\frac1{\tb}Z_{0,0}^{+,2\searrow}(2n+2;B,W)\cdot B
\label{eq:Fnexpltwo}
\end{equation}
(with some obvious notations) which, using the identities $Z_{0,0}^{+,2\searrow}(2n+2;B,W)=Z_{0,0}^+(2n+2;B,W)-Z_{0,0}^+(2n;B,W)\cdot W$ and
$B={t_\bullet}+B(B+2W)$, is easily transformed into \eqref{eq:Fnexplicit}.

\subsection{Conserved quantities for the second weighting}
We may now play a similar game for the quantities $J_n$ to obtain conserved quantities
involving the generating functions $\{P_i,Q_i\}_{i\geq 1}$. 
Let $J_n^{(d)}$ be the generating function of $\cBnd$
where each local max  (resp.\ non local max) 
receives weight $\tb$ (resp.\ $\tw$)
except for the pointed vertex that receives weight $1$. 
And let $\hat{Z}_{d,d}^+(2n;\{P_i\}_{i\geq 1},\{Q_i\}_{i\geq 1})$ be the generating
function of paths of length $2n$ starting and ending at height $d$ and staying
at height at least $d$ all along, made of elementary steps with height difference
$\pm 1$, with each descending step from height $i$ to height $i-1$ weighted
by $P_i$ if just after a descent 
and weighted by $Q_i$ if just after an ascent.

\begin{figure}
\begin{center}
\includegraphics[width=14cm]{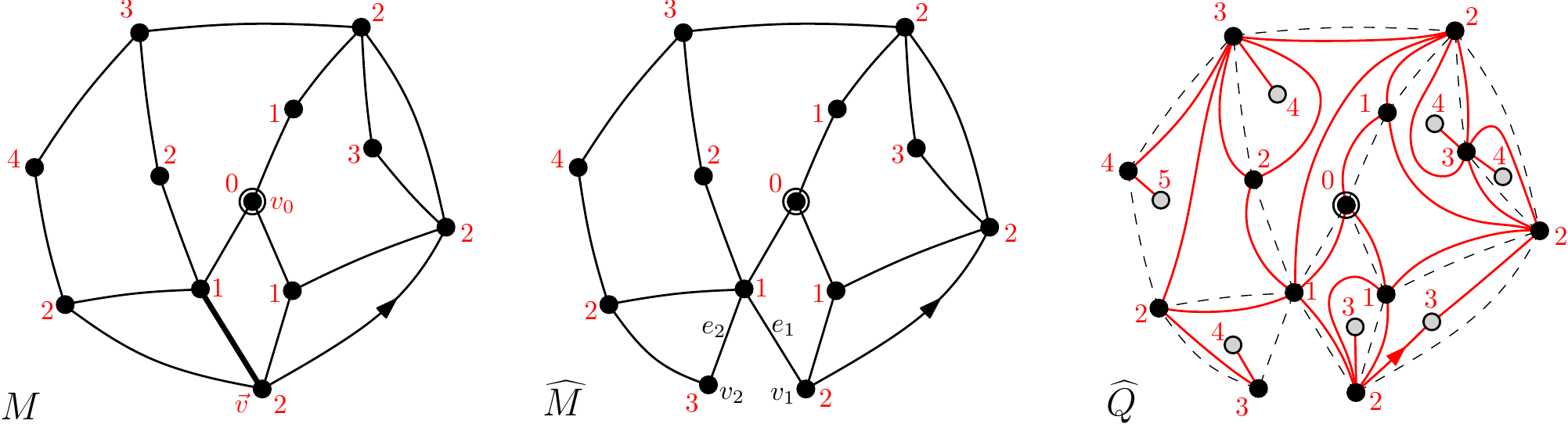}
\end{center}
\caption{Left: a pointed rooted map $M$ with a bridgeless boundary. Middle: the pointed rooted map $\hM$ obtained by cutting along the first
edge of the leftmost geodesic path from the root-vertex to the pointed vertex. Right: the associated quadrangulation $\hQ$ with
a boundary (see text).}
\label{fig:slitedgemap}
\end{figure}

Again,  
the slice decomposition described in Sect.~\ref{sec:contfracJn} for maps in $\cB_n^{(0)}$
applies more generally for maps in $\cB_n^{(d)}$ and yields
$$
J_n^{(d)}=\hat{Z}_{d,d}^+(2n;\{P_i\}_{i\geq 1},\{Q_i\}_{i\geq 1}).
$$ 
Let $\cMnd$ be the family of rooted pointed general maps with a bridgeless boundary of length $n$, 
  where the root-vertex is at distance at most $d$ from the pointed vertex $v_0$, and is at least as close from $v_0$ as any other  boundary-vertex 
(here boundary-edges are directed ccw around the map while inner edges are bi-directed; the distance-label $d(v)$ is the length of 
a shortest directed path starting from $v_0$ and ending at $v$). 
The Ambj\o rn-Budd bijection described in Sect.~\ref{sec:equality} between $\cM_n^{(0)}$ and
$\cB_n^{(0)}$ extends verbatim (using the same local rules, and having the same
pointed vertex and the same root-vertex in corresponding maps, see~\cite{AmBudd,BFG}) 
to a bijection between $\cBnd$ and $\cMnd$, so that $J_n^{(d)}$ is also the generating
function of maps in $\cMnd$ with a weight $\tb$ for each non-pointed vertex and a 
weight $\tw$ for each inner face. 

Let $\cMhnd$ be the subfamily of $\cMnd$ where the pointed vertex is different from the root-vertex, and let $\hJ_n^{(d)}$ be the generating function of the subfamily $\cBhnd$ where the weights are as in $J_n^{(d)}$. Then the Ambj\o rn-Budd bijection ensures
that $\hJ_n^{(d)}$ is also the generating function of $\cMhnd$ with a weight $\tb$
for each non-pointed vertex and a weight $\tw$ for each inner face. 
For a map $M\in\cMhnd$, let $e$ be the first edge on the leftmost geodesic path 
from the root-vertex $\vec{v}$ to the pointed vertex $v_0$ (note that all the edges on this path are inner edges). Again
we can cut along $e$ (starting from $\vec{v}$) so as to duplicate $e$ into two
edges $e_1,e_2$ (with $e_2$ before $e_1$ in ccw order around the map) and duplicate
$\vec{v}$ into two vertices $v_1,v_2$, and take $v_1$ as the new root-vertex (see Fig.~\ref{fig:slitedgemap}).
The map $\hM$ thus obtained (as opposed to the quadrangulated case we do not delete $e_2$) is a general map with a bridgeless boundary of length $n+2$.  
If we denote by $\delta_1,\ldots,\delta_{n+2}$ the distances from the pointed vertex $v_0$ 
of the successive boundary-vertices (starting with $v_1$) in ccw order around $\hM$, 
then $\delta_1$ equals the distance of $\vec{v}$ from $v_0$ in $M$ so that $\delta_1\leq d$, $\delta_i\geq \delta_1$ for $i\in\{1,\cdots,n+1\}$, 
$\delta_{i+1}\leq \delta_{i}+1$ for $i\in\{1,\cdots,n+1\}$, $\delta_{n+1}>\delta_1$ (by the effect of cutting
along the first edge of the leftmost geodesic path) 
and $\delta_{n+2}=\delta_1-1$. 
In particular if we reroot the map at the vertex between $e_1$ and $e_2$, we get a map $M'\in\cM_{n+2}^{(d-1)}$. 
We may then take the image $Q'\in\cB_{n+2}^{(d-1)}$ 
of $M'$ by the Ambj\o rn-Budd bijection, and denote by $\hQ$ the quadrangulation with boundary obtained from $Q'$ by shifting the root position by one in ccw order around $Q'$; note also that the number of local max (resp.\ non-local max) 
of $\hQ$ equals the number
of inner faces (resp.\ the number of vertices) of $M'$, which is also the number of inner faces (resp,
the number of vertices plus $1$) of $M$. 
Let again $\mathbf{d}=(d_1,\ldots,d_{2n+4})$ be the distances from the pointed vertex $v_0$ 
of the successive boundary-vertices (starting with $v_1$) in ccw order around $\hQ$. 
By the local rules of the Ambj\o rn-Budd bijection, $\mathbf{d}$ is obtained from the sequence $\delta_1,\ldots,\delta_{n+2}$ where for each $i\in\{1,\cdots,n+1\}$, we insert between $\delta_i$ and $\delta_{i+1}$ 
the subsequence (of length $\delta_i-\delta_{i+1}+1$) $\delta_{i}+1,\delta_i,\ldots,\delta_{i+1}+1$. It is then easy to check
that $\mathbf{d}$ satisfies the following conditions: $d_1=\delta_1$, $d_{i+1}=d_i\pm 1$ and $d_{i}\geq \delta_1$
for $i\in\{1,\cdots,2n+3\}$, $d_{2n+4}=\delta_1-1$, and $\mathbf{d}$ ends with $\delta_1+2,\delta_1+1,\delta_1,\delta_1-1$ (since $\delta_{n+1}>\delta_1$). 
Hence, if for $k\geq 1$ we denote by 
$\hat{Z}_{d,d}^{+,k\searrow}(2n;\{P_i\}_{i\geq 1},\{Q_i\}_{i\geq 1})$ the generating 
function defined as $\hat{Z}_{d,d}^+(2n;\{P_i\}_{i\geq 1},\{Q_i\}_{i\geq 1})$, but with 
the restriction that the $k$ last steps of the path are descending, then the slice
decomposition applied to $\widehat{Q}$ gives
$$
\hJ_n^{(d)}=\frac1{\tb}\hat{Z}_{d,d}^{+,2\searrow}(2n+2;\{P_i\}_{i\geq 1},\{Q_i\}_{i\geq 1})\cdot P_d,
$$
where the factor $\frac1{\tb}$ accounts for the  
root-vertex of $M$ being duplicated, and the factor $P_d$ accounts for the last 
descent $\delta_1,\delta_1-1$.  
Hence for each $n\geq 1$ we have the conserved quantity
\begin{equation}\label{eq:invariantJn}
J_n=\hat{Z}_{d,d}^+(2n;\{P_i\}_{i\geq 1},\{Q_i\}_{i\geq 1})-\frac1{\tb}\hat{Z}_{d,d}^{+,2\searrow}(2n+2;\{P_i\}_{i\geq 1},\{Q_i\}_{i\geq 1})\cdot P_d.
\end{equation}
Remarkably this has exactly the same form as the bicolored conserved quantities \eqref{eq:invariantFn}, 
up to changing $\{P_i,Q_i\}_{i\geq 1}$ for $\{B_i,W_i\}_{i\geq 1}$ and taking the ``hat" variants of the 
path generating functions. 
The first two invariants, $n\in\{1,2\}$, are (with $i=d-1$): for all $i\geq 1$ (with $P_0=0$)
\begin{eqnarray*}
J_1&=&Q_i-\frac1{\tb}Q_{i+1}P_iP_{i-1},\\
J_2&=&Q_i^2+Q_{i+1}P_i-\frac1{\tb}((Q_{i}+Q_{i+1})Q_{i+1}+Q_{i+2}P_{i+1})P_iP_{i-1}.
\end{eqnarray*}
As before, upon sending $d\to \infty$ in \eqref{eq:invariantJn}, we get the expression 
 \begin{equation*}
J_n=\hat{Z}_{0,0}^+(2n;P,Q)-\frac1{\tb}\hat{Z}_{0,0}^{+,2\searrow}(2n+2;P,Q)\cdot P
\end{equation*}
(with straightforward notations). Upon using $P=B$, $Q=W$ and comparing with \eqref{eq:Fnexpltwo}, this provides another
(computational) proof of the identity $J_n=F_n$ by noting that $\hat{Z}_{0,0}^+(2n;P,Q)=Z_{0,0}^+(2n;P,Q)$
and $\hat{Z}_{0,0}^{+,2\searrow}(2n+2;P,Q)=Z_{0,0}^{+,2\searrow}(2n+2;P,Q)$\footnote{
The identity $\hat{Z}_{0,0}^+(2n;P,Q)=Z_{0,0}^+(2n;P,Q)$ is easily proved by noting that the equation 
$\hat{Z}=1/(1-z (Q-P)-z P\, \hat{Z}$ which determines the generating function $\hat{Z}\equiv \sum_{n\geq 0}
\hat{Z}_{0,0}^+(2n;P,Q)z^n$ 
is identical to that,  $Z=1/(1-z Q/(1-z P\, Z))$  which determines the generating function $Z\equiv \sum_{n\geq 0}
Z_{0,0}^+(2n;P,Q)z^n$, hence $\hat{Z}=Z$. The identity $\hat{Z}_{0,0}^{+,2\searrow}(2n+2;P,Q)=Z_{0,0}^{+,2\searrow}(2n+2;P,Q)$ follows
by noting that $\hat{Z}_{0,0}^{+,2\searrow}(2n+2;P,Q)=\hat{Z}_{0,0}^+(2n+2;P,Q)-\hat{Z}_{0,0}^+(2n;P,Q)\cdot Q$
and similarly $Z_{0,0}^{+,2\searrow}(2n+2;P,Q)=Z_{0,0}^+(2n+2;P,Q)-Z_{0,0}^+(2n;P,Q)\cdot Q$.
}.
Finally, from \eqref{eq:PiQirecur}, we get $Y_1=(Q_1-P_1)=t_\circ-t_\bullet+(Q_1-P_1) Q_1$, hence
$Y_1=(t_\circ-t_\bullet)/(1-Q_1)$. Using the first conserved quantity above, we deduce $Q_1=J_1=Q-Q P^2/t_\bullet$
so that $Y_1=t_\bullet (t_\circ-t_\bullet)/(t_\bullet-t_\bullet Q+Q P^2)$ which upon expressing $t_\bullet$ and $t_\circ$
in terms of $P$ and $Q$ via \eqref{eq:BWPQ}, reproduces the expression \eqref{eq:Y1val} for $Y_1$.

 \section{Conclusion}
\label{sec:conclusion}
In this paper, we presented a comparative study of two statistical ensembles of quadrangulations. We first showed how the
corresponding slice generating functions ($B_i, W_i$ for the first ensemble and $P_i, Q_i$ for the second) appear as
coefficients of the same quantity $F(z)=J(z)$, expanded as a continued fraction in two different ways. 
The slice generating functions may then be written as bi-ratios of Hankel-type determinants and explicit
formulas may be obtained, at the price of some conjectured expression for some intermediate quantity in the second ensemble.

To conclude, we would like to emphasize that our two ensembles may be viewed, in some sense, as the two extremal
elements of a very general family of statistical ensembles as follows: by definition, the second ensemble gives 
a particular weight to those vertices which are local maxima for the distance to the root-vertex. Similarly, the first
ensemble may be viewed as the ensemble which gives a particular weight to those vertices which are local maxima 
for the \emph{distance to the root-vertex modulo $2$}. Indeed, this distance modulo $2$ is $0$ for black vertices 
(recall that the root-vertex is black) and $1$ for white vertices so that all white vertices are local maxima.
In this respect, note also that performing the passage from the quadrangulation to the general map in the bijection of Fig.~\ref{fig:quadquad2}
may be viewed as applying the Ambj\o rn-Budd rules, taking as labeling the distance modulo $2$.

Denoting by $d(v)$ the distance from a vertex $v$ to the root-vertex in a rooted quadrangulation with a boundary, we may more generally consider 
statistical ensembles which give a particular weight to those vertices which are local maxima 
for some labeling $\ell(d(v))$, with $d\mapsto \ell(d)$ some given function. Without loss of generality, 
we may set $\ell(0)=0$ and, if we wish to apply the Ambj\o rn-Budd rules to transform our quadrangulation into 
a general map, we need that $|\ell(d)-\ell(d-1)|=1$
(it also seems natural to impose that $\ell(d)$ remains non negative so that the root-vertex cannot
be a local maximum). It is likely that slice generating functions in this ensembles may appear as coefficient of $F(z)=J(z)$,
once expanded as a continued fraction with some appropriate structure, being a mixture of the Stieljes-type and
of our new encountered type. At this stage, it is interesting to notice that, in their study of \emph{finite} continued fractions
\cite{DiFKe10,DiFKe11}, Di Francesco and Kedem introduced precisely a whole family of such``mixed" fractions
as well as some passage rules on their coefficients to go from one to the other without changing the actual value of
the fraction. It is very tempting to speculate that their study may be extended to infinite continued fractions to describe our 
more general ensembles.

\appendix
 \section{A proof of the formulas \eqref{eq:YiHankel} and \eqref{eq:JJtilde}}
\begin{figure}
\begin{center}
\includegraphics[width=10cm]{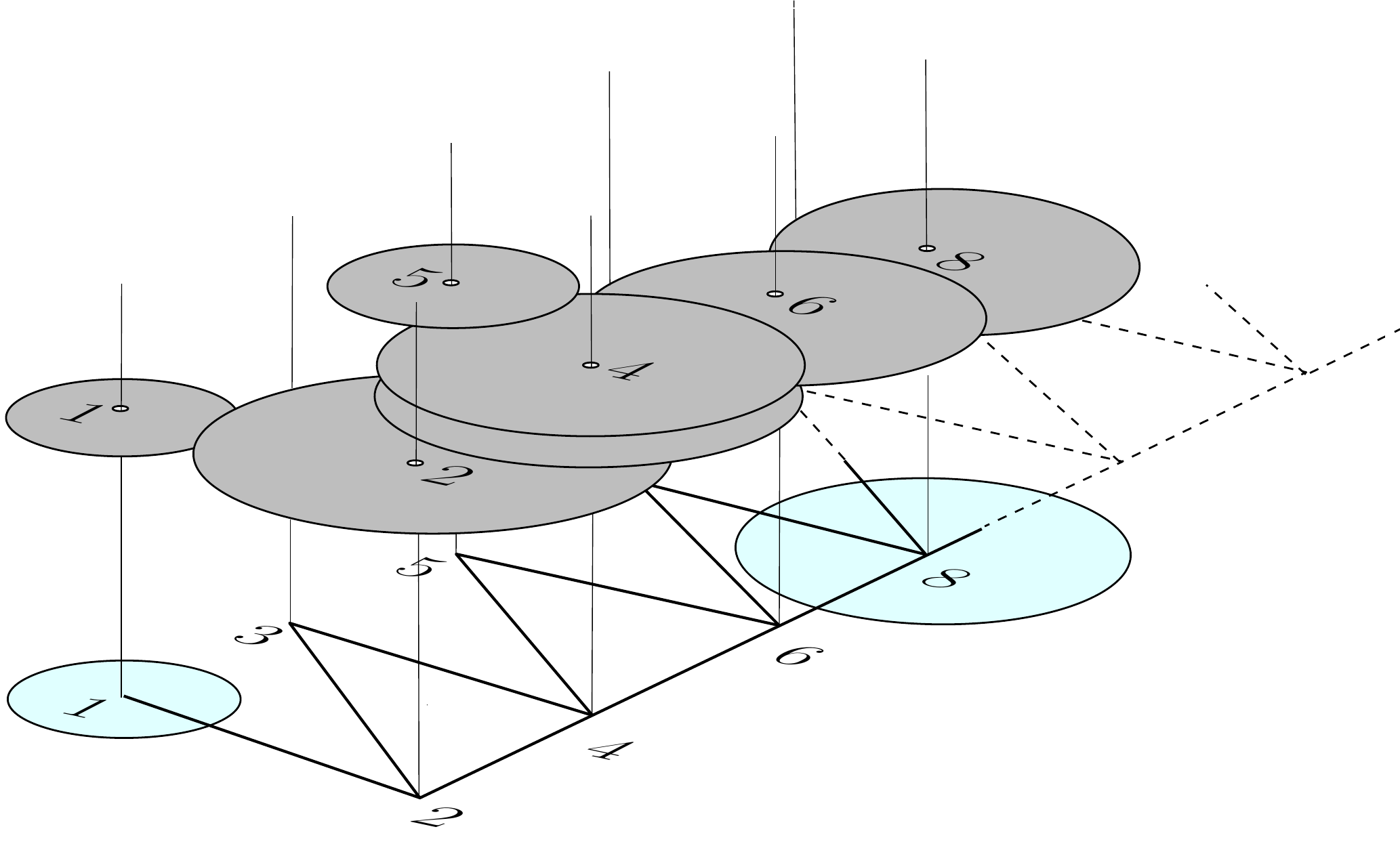}
\end{center}
\caption{A example of heap of $7$ pieces sitting on top of the graph ${\mathcal G}$ of Fig.~\ref{fig:graphG} with base $\{1,8\}$
(we indicated in light blue the ``shadow" of those pieces which can move freely and hit the vertices of the graph).
The diameter of the pieces is adjusted so that pieces sitting on top of vertices which are adjacent in ${\mathcal G}$
cannot pass through each other.}
\label{fig:heaps}
\end{figure}
As in \cite{DiFKe10,DiFKe11}, our proof of formulas  \eqref{eq:YiHankel} and \eqref{eq:JJtilde} 
 is based on the theory of \emph{heaps of pieces}. The reader is invited to consult \cite{HEAPS}
 for the basics of this theory. 
 
 Let us simply recall what we mean by a heap of pieces on a graph ${\mathcal G}$, supposedly
 connected, planar, and drawn in a \emph{horizontal} plane for simplicity. Imagine to complete the graph 
 by a set of \emph{vertical} half-lines, with a half-line starting from each vertex of the graph. Informally speaking, 
 a heap is a collection of pieces threaded along these half-lines. Each piece therefore sits on top of
 a given vertex and may move freely along the corresponding vertical half-line as long as it does not meet another piece. 
 More precisely, the pieces are supposed to be designed so that two pieces may not pass each other if they sit on top of the same 
 vertex or \emph{if they sit on top of adjacent vertices}. 
 
 Given a subset ${\mathcal B}$ of the set of vertices of ${\mathcal G}$, a heap of pieces is said to be of \emph{base} ${\mathcal B}$
 if, moving its pieces as far as possible to the bottom of the half-lines,  the set of those vertices hit by a piece 
 forms a subset of ${\mathcal B}$ (see Fig.~\ref{fig:heaps}).

 \begin{figure}
\begin{center}
\includegraphics[width=6cm]{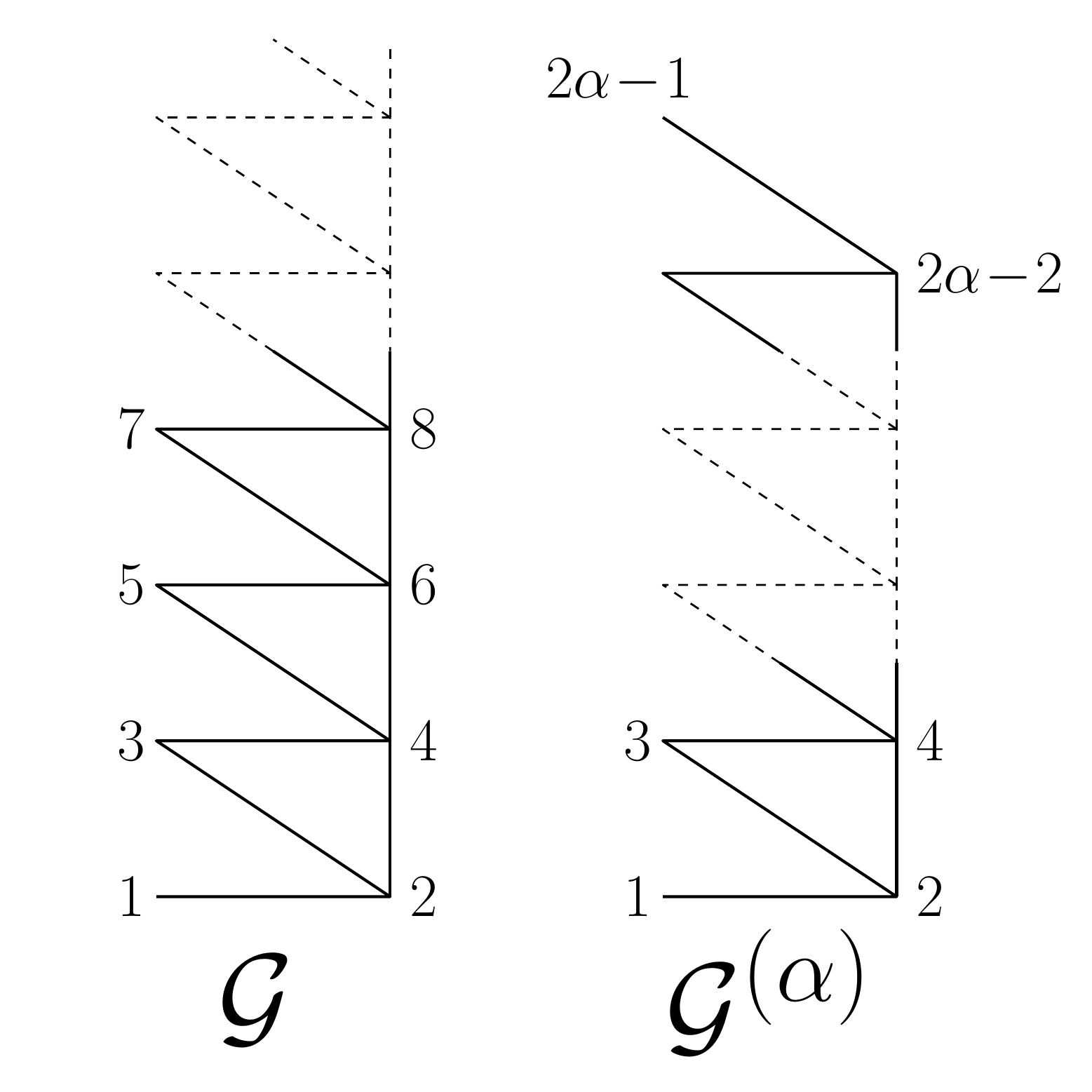}
\end{center}
\caption{The semi-infinite graph ${\mathcal G}$ and the finite graph ${\mathcal G}^{(\alpha)}$.}
\label{fig:graphG}
\end{figure}

 A fundamental remark is that, from the relation \eqref{eq:contfracnewbis}, $J(z)$ may be viewed as
 the generating function for heaps of pieces on the \emph{semi-infinite} graph ${\mathcal G}$ of 
 Fig.~\ref{fig:graphG}, with a weight $z Y_i$ per piece sitting at position $i$ along the graph, and whose base is $\{1,2\}$.
 Similarly, we may interpret $\tilde{J}(z)$ as 
 the generating function for the very same heaps, but now with a weight $z \tilde{Y}_i$ 
 per piece sitting at position $i$.
 Let us finally introduce the quantity 
 \begin{equation*}
 K(z)\equiv 1+ z\, Y_1\, J(z) = \frac{1}{\displaystyle{1-z \frac{Y_1}{\displaystyle{1-z\frac{Y_2}{\displaystyle{1-z Y_3-z \frac{Y_4}{\displaystyle{1-z Y_5-z \frac{Y_6}{\displaystyle{ 1- \cdots}}}}}}}}}}
 \end{equation*}
 which is  the generating function for heaps of pieces on the graph ${\mathcal G}$ again with a weight $z Y_i$ per piece sitting at position $i$ along the graph, but now with base $\{1\}$.

From the definition  \eqref{eq:jn} of the $j_n$'s, we have 
\begin{equation*}
K(z)=\sum_{n\geq 0} j_n z^n\ , \qquad \tilde{J}(z)=\sum_{n\geq 0} j_{-n} z^n
\end{equation*}
so that all the $j_n$'s have a direct interpretation as enumerating heap configurations made of $|n|$ pieces.

Let us now consider the analogs $J^{(\alpha)}(z)$, $K^{(\alpha)}(z)$ and $\tilde{J}^{(\alpha)}(z)$ of $J(z)$, $K(z)$ and $\tilde{J}(z)$ respectively, 
viewed as heaps generating functions, now defined on the \emph{finite} graph ${\mathcal G}^{(\alpha)}$ of Fig.~\ref{fig:graphG}. In other words, we set 
\begin{equation*}
\begin{split}
&J^{(\alpha)}(z)\equiv \frac{1}{\displaystyle{1-z Y_1-z \frac{Y_2}{\displaystyle{1-z Y_3-z \frac{Y_4}{\displaystyle{
\frac{\ddots}{\displaystyle{1-z Y_{2\alpha-3}-z \frac{Y_{2\alpha-2}}{\displaystyle{ 1- z\, Y_{2\alpha-1}}}}}}}}}}} \\
&\tilde{J}^{(\alpha)}(z)\equiv\frac{1}{\displaystyle{1-z \tilde{Y}_1-z \frac{\tilde{Y}_2}{\displaystyle{1-z \tilde{Y}_3-z \frac{\tilde{Y}_4}{\displaystyle{
\frac{\ddots}{\displaystyle{1-z \tilde{Y}_{2\alpha-3}-z \frac{\tilde{Y}_{2\alpha-2}}{\displaystyle{ 1- z\, \tilde{Y}_{2\alpha-1}}}}}}}}}}} \\
&K^{(\alpha)}(z)\equiv 1+ z\, Y_1\, J^{(\alpha)}(z)= \frac{1}{\displaystyle{1-z \frac{Y_1}{\displaystyle 1 -z \frac{Y_2}{\displaystyle{1-z Y_3-z \frac{Y_4}{\displaystyle{
\frac{\ddots}{\displaystyle{1-z Y_{2\alpha-3}-z \frac{Y_{2\alpha-2}}{\displaystyle{ 1- z\, Y_{2\alpha-1}}}}}}}}}}}}\ .\\
\end{split}
\end{equation*}
We finally define the analogs $j^{(\alpha)}_n$ of $j_n$ via
 \begin{equation*}
K^{(\alpha)}(z)=\sum_{n\geq 0} j^{(\alpha)}_n z^n\ , \qquad \tilde{J}^{(\alpha)}(z)=\sum_{n\geq 0} j^{(\alpha)}_{-n} z^n
\end{equation*}
so that $j^{(\alpha)}_n$ ($n\geq 0$) enumerates heap configurations of $n$ pieces on ${\mathcal G}^{(\alpha)}$
with weights $Y_i$ and base $\{1\}$, and $j^{(\alpha)}_{-n}$ ($n\geq 0$) enumerates heap configurations of $n$ pieces 
on ${\mathcal G}^{(\alpha)}$ with weights $\tilde{Y}_i$ and base $\{1,2\}$.

 \begin{figure}
\begin{center}
\includegraphics[width=4cm]{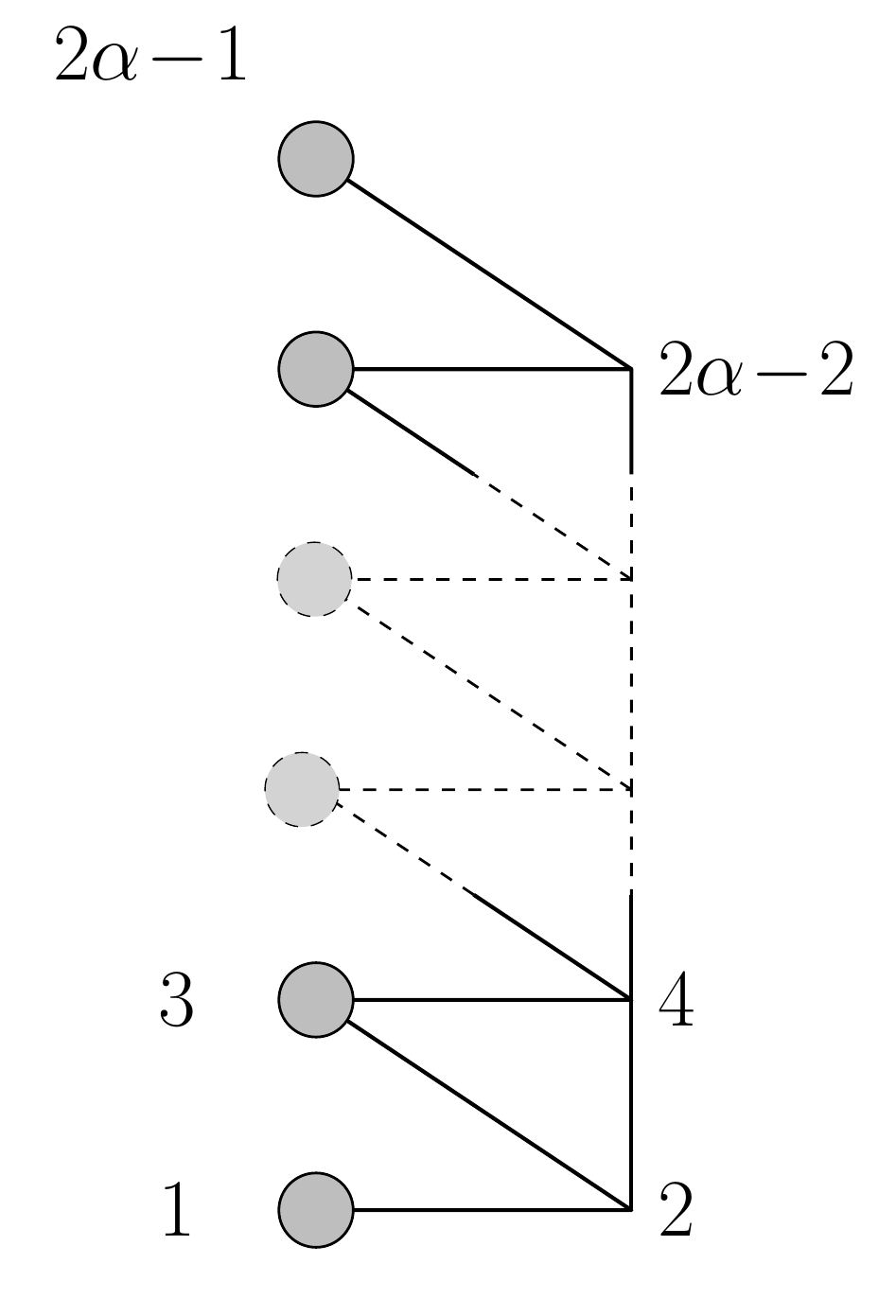}
\end{center}
\caption{The maximally occupied configuration of hard pieces on the graph ${\mathcal G}^{(\alpha)}$, made
of $\alpha$ pieces (here represented by gray circles) sitting on all the sites with an odd index.}
\label{fig:maxoccupied}
\end{figure}

It is now a standard result of the theory of heaps of pieces~\cite{HEAPS} that\footnote{
A sketch of the proof is as follows: given ${\mathcal B}$,
consider \emph{pairs} $({\mathcal H},{\mathcal HP})$ made of a heap configuration ${\mathcal H}$ of base ${\mathcal B}$ together with a configuration 
${\mathcal HP}$ of hard pieces, drawn on top of the heap. For such a pair, let $\cE$ be the set of pieces that can be moved up freely to infinity, 
and when pushed downward either are blocked by a piece (that has to be in $\mathcal H$) or hit a vertex of the base $\cB$. 
Consider the following transformation: if $\cE$ is not empty, pick the piece $p\in\cE$ of smallest index and change its status
(from $\mathcal H$ to $\mathcal HP$ if $p\in\mathcal H$, from $\mathcal HP$ to $\mathcal H$ if $p\in\mathcal HP$); if $\cE$ is empty do nothing. 
This transformation is easily seen to be an involution (which leaves $\cE$ invariant), and, if we assign a weight 
$z$ per piece in the heap and $-z$ per piece in the configuration of hard pieces, the weight is multiplied by $-1$ for each configuration
which changes under the involution. The generating function for the pairs, which is the product of the generating function for heaps with 
a weight $z$ per piece times that of configuration of hard pieces with a weight $-z$ per piece, therefore reduces to those pairs for which 
$\cE$ is empty. It is easily seen that this situation corresponds to an empty heap and a configuration of hard pieces made 
of pieces which do not belong to ${\mathcal B}$. The corresponding generating function is nothing but that of configurations of hard pieces with a 
weight $-z$ per piece not in  ${\mathcal B}$ and $0$ per piece in ${\mathcal B}$.}
\begin{equation}
\begin{split}
& K^{(\alpha)}(z)=\frac{X^{(\alpha)}(0,-z Y_2,-z Y_3,\cdots, -z Y_{2\alpha-1})}{X^{(\alpha)}(-z Y_1,-z Y_2,-z Y_3,\cdots, -z Y_{2\alpha-1})} \\
& \tilde{J}^{(\alpha)}(z)=\frac{X^{(\alpha)}(0,0,-z \tilde{Y}_3,\cdots, -z \tilde{Y}_{2\alpha-1})}{X^{(\alpha)}(-z \tilde{Y}_1,-z \tilde{Y}_2,-z \tilde{Y}_3,\cdots, -z \tilde{Y}_{2\alpha-1})} \\
\end{split}
\label{eq:heapsres}
\end{equation}
where $X^{(\alpha)}(y_1,y_2,y_3,\cdots, y_{2\alpha-1})$ denotes the generating function of \emph{hard pieces}
on the graph ${\mathcal G}^{(\alpha)}$, each piece sitting at position $i$ receiving the weight $y_i$. 
Recall that, by definition, in a configuration of hard pieces, each vertex of the graph is
occupied by at most one piece, with \emph{no two adjacent vertices occupied simultaneously}. Note that the positions of the 
$0$'s in the numerators correspond to the location of the vertices of the corresponding base of the heaps ($\{1\}$ and $\{1,2\}$ respectively).
Clearly, on the graph ${\mathcal G}^{(\alpha)}$, we can put at most $\alpha$ hard pieces. Moreover, this maximal situation is achieved
by a single configuration with all sites with odd index occupied (see Fig.~\eqref{fig:maxoccupied}). The quantity 
$X^{(\alpha)}(-z Y_1,-z Y_2,-z Y_3,\cdots, -z Y_{2\alpha-1})$ is therefore a polynomial of degree $\alpha$ in $z$ 
that we write
\begin{equation*}
X^{(\alpha)}(-z Y_1,-z Y_2,-z Y_3,\cdots, -z Y_{2\alpha-1})=\sum_{m=0}^\alpha (-z)^m X^{(\alpha)}_m(Y_1,Y_2,Y_3,\cdots,Y_{2\alpha-1})
\end{equation*}
where $X^{(\alpha)}_m(y_1,y_2,y_3,\cdots, y_{2\alpha-1})$ denotes the generating function of exactly $m$ \emph{hard pieces}
on the graph ${\mathcal G}^{(\alpha)}$, each piece sitting at position $i$ receiving the weight $y_i$.
Clearly, both $X^{(\alpha)}(0,-z Y_2,-z Y_3,\cdots, -z Y_{2\alpha-1})$ and
$X^{(\alpha)}(0,0,-z \tilde{Y}_3,\cdots, -z \tilde{Y}_{2\alpha-1})$ are polynomials of degree $\alpha-1$ in $z$.

 \begin{figure}
\begin{center}
\includegraphics[width=6cm]{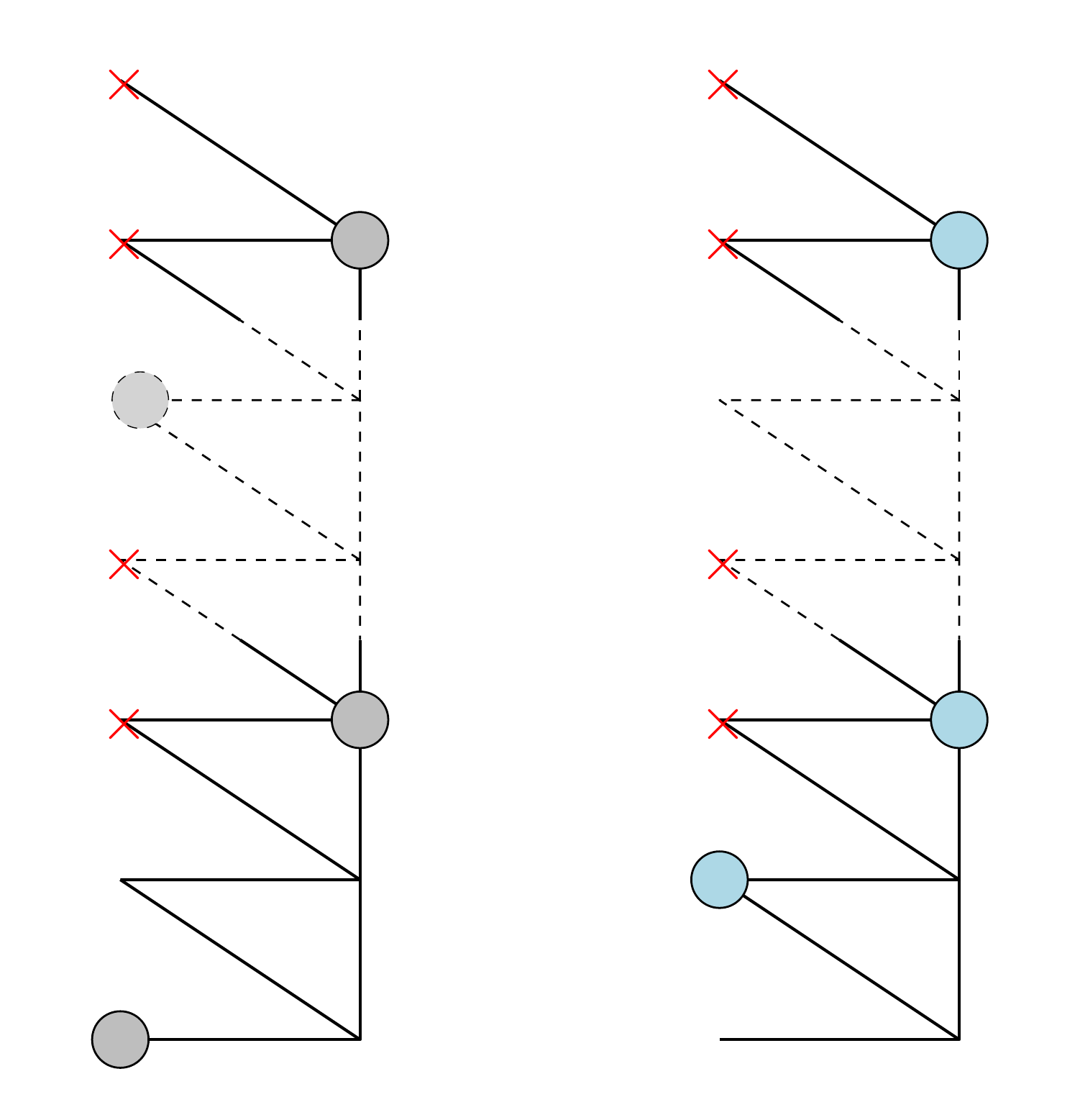}
\end{center}
\caption{Left: an example of configuration $C$ of hard pieces (represented in gray) on ${\mathcal G}^{(\alpha)}$.
Right: the associated configuration $\tilde{C}$ of hard pieces (represented in light blue), obtained by keeping
the particles sitting on even vertices and, in the ensemble of odd vertices \emph{which are not adjacent to the occupied 
even vertices}, exchanging the occupied and un-occupied sites.}
\label{fig:YYtilde}
\end{figure}

Let us now come to our fundamental identities. We have
\begin{equation}
\begin{split}
&  X^{(\alpha)}_m =X^{(\alpha)}_\alpha\ \tilde{X}^{(\alpha)}_{\alpha-m}\\
& X^{(\alpha)}_m(0) =X^{(\alpha)}_\alpha\ \Big(\tilde{X}^{(\alpha)}_{\alpha-m}
-\tilde{X}^{(\alpha)}_{\alpha-m}(0,0)\Big)\\
\end{split}
\label{eq:XXt}
\end{equation}
with the short-hand notations
\begin{equation*}
\begin{split}
& X^{(\alpha)}_m\equiv X^{(\alpha)}_m(Y_1,Y_2,Y_3,\cdots,Y_{2\alpha-1}) \\
& X^{(\alpha)}_m(0)\equiv X^{(\alpha)}_m(0,Y_2,Y_3,\cdots,Y_{2\alpha-1}) \\
& \tilde{X}^{(\alpha)}_{m}\equiv X^{(\alpha)}_m(\tilde{Y_1},\tilde{Y_2},\tilde{Y}_3,\cdots, \tilde{Y}_{2\alpha-1})\\
& \tilde{X}^{(\alpha)}_{m}(0,0)\equiv X^{(\alpha)}_m(0,0,\tilde{Y}_3,\cdots, \tilde{Y}_{2\alpha-1})\ .\\
\end{split}
\end{equation*}
To explain these identities, let us analyze the structure of a configuration $C$ of hard pieces on ${\mathcal G}^{(\alpha)}$.
In $C$, a number $k$ of pieces occupy even sites $2j_1,2j_2,\cdots,2j_k$ with $1\leq j_1<j_2<\cdots<j_k\leq \alpha-1$ and $j_{\ell+1}-j_\ell>1$
for $\ell=1,\cdots, k-1$. The set of available odd sites is ${\mathcal Odd}=\{1,3,5,\cdots,2\alpha-1\}\setminus \{2j_1-1,2j_1+1,2j_2-1,2j_2+1,\cdots,
2j_k-1,2j_k+1\}$ and satisfies $|{\mathcal Odd}|=\alpha-2k$. A number $k'$ of pieces occupy a subset $\{2j'_1-1,2j'_2-1,\cdots,2j'_{k'}-1\}$ of this set. In $X^{(\alpha)}_m$
(corresponding to a situation where $k+k'=m$), any occupied site $i$ receives the weight $Y_i$, so the weight of the configuration $C$ is
$Y_{2j_1}Y_{2j_2}\cdots Y_{2j_k}\times Y_{2j'_1-1}Y_{2j'_2-1}\cdots Y_{2j'_{k'}-1}$. Let us now consider instead the configuration $\tilde{C}$ where, again, the sites 
$2j_1,2j_2,\cdots,2j_k$ are occupied but now the complementary of $\{2j'_1-1,2j'_2-1,\cdots,2j'_{k'}-1\}$ in ${\mathcal Odd}$ (namely ${\mathcal Odd}\setminus \{2j'_1-1,2j'_2-1,\cdots,2j'_{k'}-1\}$) is covered by pieces. Clearly, going in from $C$ to $\tilde{C}$ provides a bijection between 
configurations $C$ with $k+k'=m$ pieces and configurations $\tilde{C}$ with $k+(\alpha-2k)-k'=\alpha-m$ pieces. In $\tilde{X}_{\alpha-m}$,
the configuration $\tilde{C}$ receives the weight 
\begin{equation*}
\begin{split}
\hskip -1.cm \tilde{Y}_{2j_1}\tilde{Y}_{2j_2}\cdots \tilde{Y}_{2j_k} \times \frac{\prod\limits_{i\in {\mathcal Odd}}\tilde{Y}_i}{\tilde{Y}_{2j'_1-1}\tilde{Y}_{2j'_2-1}\cdots \tilde{Y}_{2j'_{k'}-1}}
&= \frac{Y_{2j_1}Y_{2j_2}\cdots Y_{2j_k}}{Y_{2j_1-1}Y_{2j_1+1}Y_{2j_2-1}Y_{2j_2+1}\cdots Y_{2j_k-1}Y_{2j_k+1}}
\\ & \qquad \qquad \qquad \times \frac{Y_{2j'_1-1}Y_{2j'_2-1}\cdots Y_{2j'_{k'}-1}}{\prod\limits_{i\in {\mathcal Odd}}Y_i}\\
&=\frac{Y_{2j_1}Y_{2j_2}\cdots Y_{2j_k}\times Y_{2j'_1-1}Y_{2j'_2-1}\cdots Y_{2j'_{k'}-1}}{Y_1Y_3Y_5\cdots Y_{2\alpha-1}}\\
&=\frac{Y_{2j_1}Y_{2j_2}\cdots Y_{2j_k}\times Y_{2j'_1-1}Y_{2j'_2-1}\cdots Y_{2j'_{k'}-1}}{X^{(\alpha)}_\alpha}\\
\end{split}
\end{equation*}
since $X^{(\alpha)}_\alpha=Y_1Y_3Y_5\cdots Y_{2\alpha-1}$.
From the bijection $C\mapsto \tilde{C}$, we therefore deduce immediately the first equality in \eqref{eq:XXt}. To get the second equality,
we note that $X^{(\alpha)}_m(0)$ enumerates configurations $C$ with $m$ pieces such that \emph{the site $1$ is not occupied by a piece}. Two situations may
then occur: either site $2$ is occupied or not. In the first case, the
bijection $C\mapsto \tilde{C}$ will generate a configuration $\tilde{C}$ where site $2$ is occupied (and site $1$ does not belong to 
${\mathcal Odd}$) while in the second case, 
it will generate a configuration where site $2$ is empty and site $1$ (which belongs to ${\mathcal Odd}$) is necessarily 
occupied (since it was empty in $C$ and the empty and occupied sites get exchanged in the bijection for those
odd sites belonging to ${\mathcal Odd}$). To summarize, in the configuration $\tilde{C}$, either site $1$ or site $2$ must
be occupied. The restriction of $\tilde{X}^{(\alpha)}_{\alpha-m}$ to these configurations yields $\tilde{X}^{(\alpha)}_{\alpha-m}-\tilde{X}^{(\alpha)}_{\alpha-m}(0,0)$,
hence the second equality.

From \eqref{eq:XXt}, we deduce
\begin{equation*}
\begin{split}
&\hskip -1.3cm  X^{(\alpha)}\left(-\frac{Y_1}{z},-\frac{ Y_2}{z},-\frac{Y_3}{z},\cdots, -\frac{Y_{2\alpha-1}}{z}\right)\!=\!\left(-\frac{1}{z}\right)^{\alpha} \!\!X^{(\alpha)}_\alpha\!\cdot\! 
 X^{(\alpha)}(-z \tilde{Y}_1,-z \tilde{Y}_2,-z \tilde{Y}_3,\cdots, -z \tilde{Y}_{2\alpha-1})\\
&\hskip -1.3cm  X^{(\alpha)}\left(0,-\frac{ Y_2}{z},-\frac{ Y_3}{z},\cdots, -\frac{ Y_{2\alpha-1}}{z}\right)\!=\!\left(-\frac{1}{z}\right)^{\alpha}\!\! X^{(\alpha)}_\alpha\!\cdot\!
\Big(X^{(\alpha)}(-z \tilde{Y}_1,-z \tilde{Y}_2,-z \tilde{Y}_3,\cdots, -z \tilde{Y}_{2\alpha-1})\\
& \hskip 7.7cm -
X^{(\alpha)}(0,0,-z \tilde{Y}_3,\cdots, -z \tilde{Y}_{2\alpha-1})\Big)\\
\end{split}
\end{equation*}
and therefore, taking the ratio of the two lines and using \eqref{eq:heapsres},
\begin{equation*}
1+\frac{Y_1}{z}J^{(\alpha)}\left(\frac{1}{z}\right)=K^{(\alpha)}\left(\frac{1}{z}\right)=1-\tilde{J}^{(\alpha)}(z)\ .
\end{equation*}
The finite continued fraction case of Sect.~\ref{sec:finite} corresponds precisely to a situation
where $J(z)=J^{(\alpha)}(z)$ and $\tilde{J}(z)=\tilde{J}^{(\alpha)}(z)$. The above formula explains the 
first identity in \eqref{eq:JJtilde} while the second identity is guaranteed by the relation $\tilde{J}^{(\alpha)}(z)\to 1$ when $z\to 0$.
This concludes the proof of \eqref{eq:JJtilde}.

We now prove 
\eqref{eq:YiHankel} by computing explicitly the determinants $H_i^{(0)}$ and
$H_i^{(1)}$ in terms of the $Y_i$'s. More precisely, let us show that, for $i\geq 1$,
\begin{equation}
\begin{split}
H_i^{(0)}&\equiv \left\vert 
\begin{matrix}
j_{-(i-1)} & \cdots & \cdots & j_0 \\
\vdots & \rddots & \rddots  & j_1\\
\vdots & \rddots & \rddots & \vdots \\
j_0 & j_1 & \cdots & j_{i-1} \\
\end{matrix}
\right\vert = \left(\frac{Y_2}{Y_3}\right)^{i-1}\left(\frac{Y_4}{Y_5}\right)^{i-2}\cdots\left(\frac{Y_{2i-4}}{Y_{2i-3}}\right)^{2}\left(\frac{Y_{2i-2}}{Y_{2i-1}}\right) \\
& \\
H_i^{(1)}&\equiv \left\vert 
\begin{matrix}
j_{-(i-2)} & \cdots & \cdots & j_1 \\
\vdots & \rddots & \rddots  & j_2\\
\vdots & \rddots & \rddots & \vdots \\
j_1 & j_2 & \cdots & j_{i} \\
\end{matrix}
\right\vert =  Y_1\, Y_3\,  Y_5 \cdots Y_{2i-1}\, H_i^{(0)}\ .\\
\end{split}
\label{eq:HHexpl}
\end{equation}
Once these formulas are proved, the relations \eqref{eq:YiHankel} indeed follow immediately.

A first crucial point is the existence of a linear relation between the $j^{(\alpha)}_n$'s, namely
\begin{equation}
\sum_{m=0}^\alpha (-1)^m X^{(\alpha)}_m\,  j^{(\alpha)}_{n-m}=0\ \hbox{for all integers}\ n\  .
\label{eq:linearrel}
\end{equation} 
Indeed, writing the first identity in \eqref{eq:heapsres} as $K^{(\alpha)}(z)X^{(\alpha)}(-z Y_1,-z Y_2,-z Y_3,\cdots, -z Y_{2\alpha-1})=X^{(\alpha)}(0,-z Y_2,-z Y_3,\cdots, -z Y_{2\alpha-1})$ and extracting the term of order $z^n$, we immediately see that \eqref{eq:linearrel} holds
for any positive integer $n\geq \alpha$ since $X^{(\alpha)}(0,-z Y_2,-z Y_3,\cdots, -z Y_{2\alpha-1})$ is a polynomial of degree $\alpha-1$.
Similarly, writing the second identity in \eqref{eq:heapsres} as $\tilde{J}^{(\alpha)}(z) \times X^{(\alpha)}(-z \tilde{Y}_1,-z \tilde{Y}_2,\cdots, -z \tilde{Y}_{2\alpha-1})=X^{(\alpha)}(0,0,-z \tilde{Y}_3,\cdots, -z \tilde{Y}_{2\alpha-1})$, a polynomial of degree $\alpha-1$, we find that 
\begin{equation*}
\sum_{m=0}^\alpha (-1)^m \tilde{X}^{(\alpha)}_m\,  j^{(\alpha)}_{-(n'-m)}=0= \frac{(-1)^\alpha}{X^{(\alpha)}_\alpha}
\sum_{m'=0}^\alpha (-1)^{m'} X^{(\alpha)}_{m'}\,  j^{(\alpha)}_{(\alpha-n')-m'}
\end{equation*} 
for $n'\geq \alpha$. Here we have set $m'=\alpha-m$ and used $\tilde{X}^{(\alpha)}_{\alpha-m'}=X^{(\alpha)}_{m'}/X^{(\alpha)}_\alpha$.
Setting $n=\alpha-n'\leq 0$, we deduce that \eqref{eq:linearrel} also holds for any non-positive integer. It remains to show that 
it is valid in the range $1\leq n \leq \alpha-1$.  For $n$ in this range, we have
\begin{equation*}
\begin{split}
\sum_{m=0}^\alpha (-1)^m X^{(\alpha)}_m\,  j^{(\alpha)}_{n-m}&=
\sum_{m=0}^n (-1)^m X^{(\alpha)}_m\,  j^{(\alpha)}_{n-m}+
\sum_{m=n+1}^\alpha (-1)^m X^{(\alpha)}_m\,  j^{(\alpha)}_{n-m} \\
& = (-1)^n X^{\alpha}_n(0)+(-1)^\alpha\sum_{p=0}^{\alpha-n-1} (-1)^p X^{(\alpha)}_{\alpha-p}\,  j^{(\alpha)}_{p-(\alpha-n)}\\
& = (-1)^n X^{\alpha}_n(0)+(-1)^\alpha X_\alpha^{(\alpha)} \sum_{p=0}^{\alpha-n-1} (-1)^p \tilde{X}^{(\alpha)}_{p}\,  j^{(\alpha)}_{-((\alpha-n)-p)}\ .\\
\end{split}
\end{equation*}
Now since $1\leq \alpha-n\leq \alpha-1$, we also have
\begin{equation*}
\sum_{p=0}^{\alpha-n} (-1)^p \tilde{X}^{(\alpha)}_p\,  j^{(\alpha)}_{-((\alpha-n)-p)}=(-1)^{\alpha-n} \tilde{X}^{(\alpha)}_{\alpha-n}(0,0)
\end{equation*}
so that we eventually get
\begin{equation*}
\begin{split}
\hskip -1.cm \sum_{m=0}^\alpha (-1)^m X^{(\alpha)}_m\,  j^{(\alpha)}_{n-m}&=(-1)^n X^{\alpha}_n(0)
+(-1)^\alpha X_\alpha^{(\alpha)} \left((-1)^{\alpha-n} \tilde{X}^{(\alpha)}_{\alpha-n}(0,0)-(-1)^{\alpha-n} \tilde{X}^{(\alpha)}_{\alpha-n} \right)\\
&=  (-1)^n\left( X^{\alpha}_n(0)-X_\alpha^{(\alpha)} ( \tilde{X}^{(\alpha)}_{\alpha-n}- \tilde{X}^{(\alpha)}_{\alpha-n}(0,0))\right)=0\ .\\
\end{split}
\end{equation*}
The linear relation \eqref{eq:linearrel} therefore holds for all integers $n$, as stated.

Let us now come to the computation of $H_i^{(0)}$ and $H_i^{(1)}$. Since $j_n$, $n\geq 1$, enumerates heaps of $n$ pieces 
on the graph ${\mathcal G}$ with base $\{1\}$, the pieces cannot reach sites with index more than $1$ for $n=1$ and $2n-2$ for 
$n\geq 2$. In other words, $j_n$ enumerate heaps which ``live" on ${\mathcal G}^{(n)}$, therefore on ${\mathcal G}^{(i-1)}$
for all $n\leq i-1$. As for $j_{n}$, $n\leq - 1$, it enumerates heaps of $|n|$ pieces 
on the graph ${\mathcal G}$ with base $\{1,2\}$ so that the pieces cannot reach sites with index more than $2|n|$,
therefore ``live" on ${\mathcal G}^{(|n|+1)}$, therefore on ${\mathcal G}^{(i-1)}$
for all $|n|\leq i-2$. In other words, we have
\begin{equation*}
\begin{split}
& j_n=j_n^{(i-1)}\ \hbox{for}\ n=0,1,2,\cdots, i-1\\
& j_{-n}=j_{-n}^{(i-1)}\ \hbox{for}\ n=0,1,2,\cdots, i-2\\
\end{split}
\end{equation*}
In the determinant $H_i^{(0)}$, the only term which does not ``live" on ${\mathcal G}^{(i-1)}$ is $j_{-(i-1)}$
and it is easily seen that
\begin{equation*}
j_{-(i-1)}=j^{(i-1)}_{-(i-1)}+ \tilde{Y}_2\tilde{Y_4}\cdots\tilde{Y}_{2i-2}
\end{equation*}
with an additional term corresponding to the unique heap that hits position $2i-2$.
Using the linear relation \eqref{eq:linearrel} for $\alpha=i-1$, we may thus rewrite $H_i^{(0)}$ as
\begin{equation*}
H_i^{(0)}= \left\vert 
\begin{matrix}
 \tilde{Y}_2\tilde{Y_4}\cdots\tilde{Y}_{2i-2} & j_{-(i-2)} & \cdots & j_0 \\
0 & j_{-(i-3)} & \rddots  & j_1\\
\vdots & \rddots & \rddots & \vdots \\
0 & j_1 & \cdots & j_{i-1} \\
\end{matrix} \right\vert
= \tilde{Y}_2\tilde{Y_4}\cdots\tilde{Y}_{2i-2}\, H_{i-1}^{(1)}\ .
\end{equation*}
Similarly, in the determinant $H_i^{(1)}$, the only term which does not ``live" on ${\mathcal G}^{(i-1)}$ is $j_{i}$
and it is easily seen that
\begin{equation*}
j_{i}=j^{(i-1)}_{i}+ Y_1(Y_2 Y_4\cdots Y_{2i-2})
\end{equation*}
with again an additional term corresponding to the unique heap that hits position $2i-2$.
We may thus rewrite $H_i^{(1)}$ as
\begin{equation*}
H_i^{(1)}\equiv \left\vert 
\begin{matrix}
j_{-(i-2)} & \cdots & j_0 & 0 \\
\vdots & \rddots & \rddots  &\vdots\\
j_0 & \rddots & \rddots &0 \\
j_1 & \cdots & j_{i-1} & Y_1(Y_2 Y_4\cdots Y_{2i-2}) \\
\end{matrix}
\right\vert = Y_1(Y_2 Y_4\cdots Y_{2i-2})\, H_{i-1}^{(0)}\ .
\end{equation*}
Combining the two above formulas and replacing the $\tilde{Y}_i$'s by their value
in terms of the $Y_i$'s, we deduce the recursion relation 
\begin{equation*}
H_i^{(0)}=\left(\frac{Y_2}{Y_3}\frac{Y_4}{Y_5}\cdots\frac{Y_{2i-4}}{Y_{2i-3}}\right)^2 \frac{Y_{2i-2}}{Y_{2i-1}}\ H^{(0)}_{i-2}
\end{equation*}
for $i\geq 3$ with initial conditions $H^{(0)}_1=1$ and $H^{(0)}_2=j_{-1}j_{1}-j_0^2=(\tilde{Y}_1+\tilde{Y_2}) Y_1-1= (Y_2/Y_3)$.
The first line of eq~\eqref{eq:HHexpl} follows immediately. As for the second line, it follows from
\begin{equation*}
\hskip -1.2cm\frac{H^{(1)}_i}{H^{(0)}_i}=Y_1(Y_2Y_4\cdots Y_{2i-2})\frac{H^{(0)}_{i-1}}{H^{(0)}_i}= Y_1(Y_2Y_4\cdots Y_{2i-2})
\left(\frac{Y_3}{Y_2}\frac{Y_5}{Y_4}\cdots\frac{Y_{2i-1}}{Y_{2i-2}}\right)=Y_1 Y_3 Y_5 \cdots Y_{2i-1}\ .
\end{equation*}
The above derivation of Eq.~\eqref{eq:YiHankel} extends verbatim to the case of the finite continued fraction
of Sect.~\ref{sec:finite} by limiting to $i \leq \alpha-1$ the range of allowed values for the index  $i$ in $H_i^{(0)}$ and $H_i^{(1)}$.
This range is precisely what is needed to compute $Y_1,Y_2,\cdots,Y_{2\alpha-1}$.

\section{A proof of the formulas \eqref{eq:linearA}--\eqref{eq:linearD}}
The quantities $Z(z;P,Q)$, $\tilde{Z}(z;P,Q)$ are specializations of $J(z)$ and $\tilde{J}(z)$ (viewed as defined
from the $Y_i$'s through their continued fraction expansions) to the case where
\begin{equation*}
Y_{2i-1}=Y\ , \qquad Y_{2i}=P\ , \qquad \tilde{Y}_{2i-1}=\frac{1}{Y}=\tilde{Y}$\ , \qquad $\tilde{Y}_{2i}=\frac{\displaystyle{P}}{\displaystyle{Y^2}}=\tilde{P}
\end{equation*}
for all $i\geq 1$.
Consequently, $Z_n$, $\tilde{Z}_n$ and $k_n$ are the corresponding specializations of $J_n$, $\tilde{J}_n$ and $j_n$. 
The analysis of Appendix A applies to arbitrary $Y_i$'s. In particular, $k_n$ has a direct interpretation
in term of heaps of pieces on the graph ${\mathcal G}$ for all $n$.
For $n>0$, $k_n$ enumerates heaps of $n$ pieces of base $\{1\}$, with weights $Y$ and $P$ for pieces on
odd or even sites respectively. 
For $n<0$, $k_n$ enumerates heaps of $|n|$ pieces of base $\{1,2\}$, with weights $\tilde{Y}$ and $\tilde{P}$ for pieces on
odd or even sites respectively. 
Let us thus introduce the quantities $k_n^{(\alpha)}$ (analogs of $j_n^{(\alpha)}$) corresponding to
a restriction of the heaps to the graph   ${\mathcal G}^{(\alpha)}$.
From \eqref{eq:linearrel}, we deduce immediately
\begin{equation*}
\sum_{m=0}^\alpha (-1)^m x^{(\alpha)}_m\,  k^{(\alpha)}_{n-m}=0\ \hbox{for all integers}\ n\  .
\end{equation*} 
where $x^{(\alpha)}_m=X^{(\alpha)}_m(Y,P,Y,P,\cdots, Y)$ is the generating function of configurations of exactly $m$ hard pieces
on the graph ${\mathcal G}^{(\alpha)}$. Setting $\alpha=i-1$, this equation reads equivalently
\begin{equation*}
\sum_{m=0}^{i-1} k^{(i-1)}_{n-i+m} (-1)^m  \frac{x^{(i-1)}_{i-1-m}}{x^{(i-1)}_{i-1}}=0\ \hbox{for all integers}\ n\  .
\end{equation*}
Now, from their heap interpretation, it is clear that $k_n=k_n^{(i-1)}$ for $n=0,1,\cdots, i-1$ as well as for 
$n=-1,-2,\cdots, -(i-2)$. For $2\leq n\leq i$, all the $k^{(i-1)}_p$'s appearing in the above formula may
thus be
replaced by $k_p$'s and \eqref{eq:linearA} follows. For $n=1$, the only term which gets out of the 
graph ${\mathcal G}^{(i-1)}$ is for $m=0$, since $k^{(i-1)}_{-(i-1)}=k_{-(i-1)}-\tilde{P}^{i-1}$ (the two indeed differ
by the contribution of the heap made of one piece on each even site from $2$ to $2(i-1)$). This
explains the right hand side $\tilde{P}^{i-1}=P^{i-1}/Y^{2(i-1)}$ in the first line of \eqref{eq:linearB}.
For $n=i+1$, the only term which gets out of the 
graph ${\mathcal G}^{(i-1)}$ is for $m=i-1$, since $k^{(i-1)}_{i}=k_{i} - Y P^{i-1}$ (the two indeed differ
by the contribution of the heap made of one piece on site $1$ and one piece on each even site from $2$ to $2(i-1)$). This
explains the right hand side $(-1)^{i-1}Y P^{i-1}(x^{(i-1)}_{0}/x^{(i-1)}_{i-1})=(-1)^{i-1}P^{i-1}/Y^{i-2}$ 
in the second line of \eqref{eq:linearB}. 

 \begin{figure}
\begin{center}
\includegraphics[width=3cm]{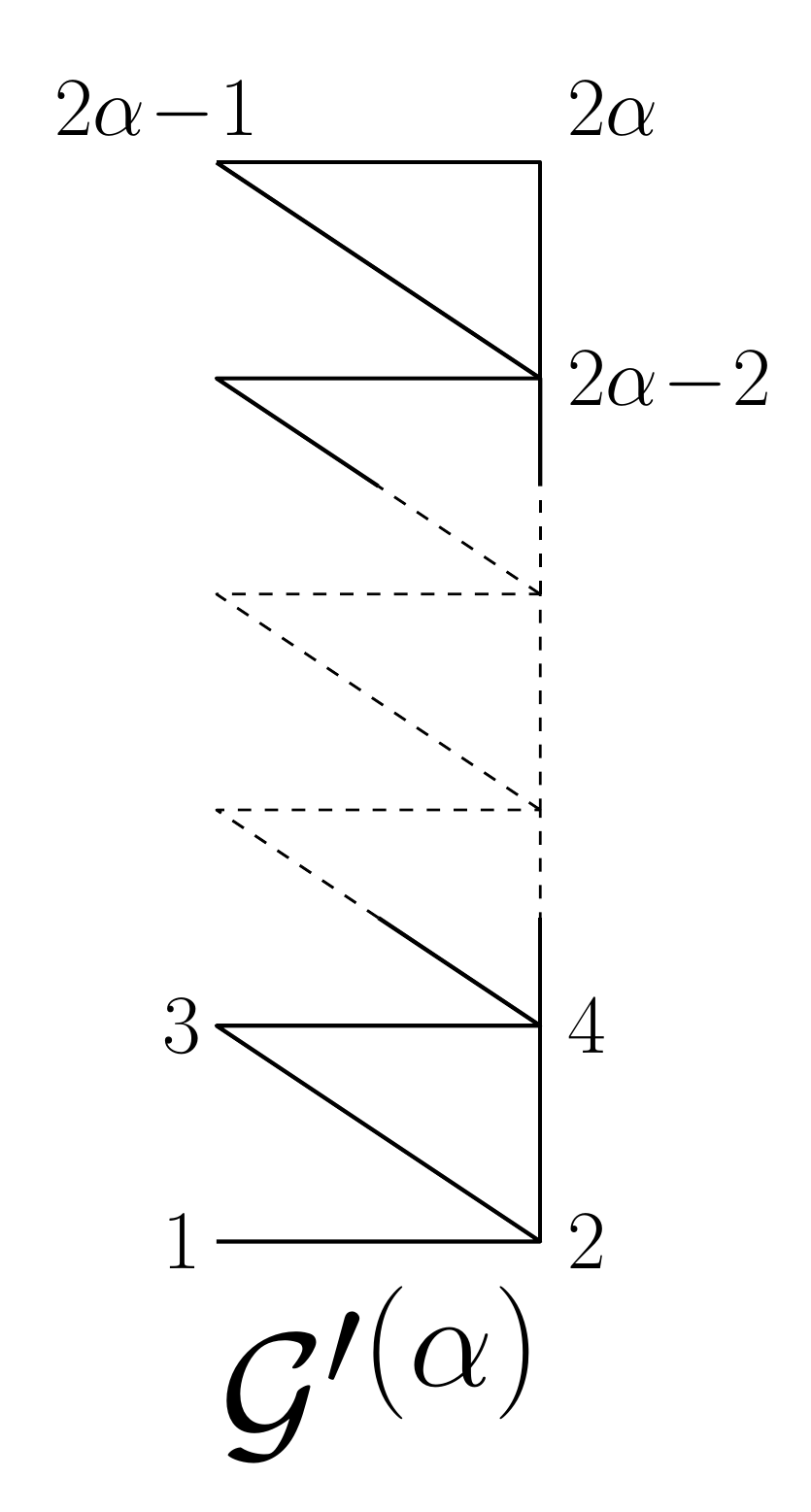}
\end{center}
\caption{The graph ${\mathcal G}'^{(\alpha)}$.}
\label{fig:graphGprime}
\end{figure}

Eqs.~\eqref{eq:linearC} and \eqref{eq:linearD} can be proved in the same way but their proof now relies
on a restriction of the heaps to the graph ${\mathcal G}'^{(\alpha)}$ of Fig.~\ref{fig:graphGprime}.
Let us first analyze the heap generating functions on this graph: denoting, for $n>0$, the generating function $j'^{(\alpha)}_n$ 
of heaps of $n$ pieces on ${\mathcal G}'^{(\alpha)}$, of base $\{1\}$ and with weight $Y_i$ per piece sitting on site $i$.
It is easily seen that $j'^{(\alpha)}_n$ also corresponds to enumerating heaps of $n$ pieces 
on the graph ${\mathcal G}^{(\alpha)}$ provided we assign, instead of $Y_{2\alpha-1}$, a weight 
$Y'_{2\alpha-1}\equiv Y_{2\alpha-1}+Y_{2\alpha}$ to pieces sitting at position $2\alpha-1$. The same remark holds for
configurations of $m$ hard pieces on ${\mathcal G}'^{(\alpha)}$ which are enumerated by $X^{(\alpha)}_m(Y_1,Y_2,\cdots, Y_{2\alpha-2},
Y'_{2\alpha-1})$. In order to use directly our previous results (obtained for ${\mathcal G}^{(\alpha)}$), we are thus lead to define $j'^{(\alpha)}_{n}$ for $n<0$ as enumerating heaps of $n$ pieces 
with base $\{1,2\}$ on the graph ${\mathcal G}^{(\alpha)}$, with weights $\tilde{Y}_i$ built via the same expression \eqref{eq:tildeYi} as before 
with $Y_{2\alpha-1}$ replaced by $Y'_{2\alpha-1}$.  Getting back to ${\mathcal G}'^{(\alpha)}$,  $j'^{(\alpha)}_{n}$ for $n<0$ therefore enumerates heaps of $n$ pieces 
with base $\{1,2\}$ on this graph, with weights $ \tilde{Y_i}$ as in \eqref{eq:tildeYi} for pieces on sites $1,2,3,\cdots, 2\alpha-3$ and the special weights $Y_{2\alpha-2}/(Y_{2\alpha-3}(Y_{2\alpha-1}+Y_{2\alpha}))$ for pieces on the site $2\alpha-2$, 
$1/(Y_{2\alpha-1}+Y_{2\alpha})$ for pieces on the site $2\alpha-1$ and $0$ for pieces on the site $2\alpha$.
With this definition, we have the analog of \eqref{eq:linearrel}, namely
\begin{equation*}
\sum_{m=0}^\alpha (-1)^m X'^{(\alpha)}_m\,  j'^{(\alpha)}_{n-m}=0\ \hbox{for all integers}\ n\  .
\end{equation*} 
where $X'^{(\alpha)}_m$ enumerates configurations of $m$ hard pieces on ${\mathcal G}'^{(\alpha)}$.

Let us now specialize this result upon introducing, for $n>0$, the generating function $k'^{(\alpha)}_n$ of heaps of $n$ pieces on ${\mathcal G}'^{(\alpha)}$, of base $\{1\}$
and with weights $Y$ and $P$ for pieces on odd or even sites respectively. From the above discussion, 
for $n<0$, $k'^{(\alpha)}_{-n}$ must be defined as enumerating heaps of $n$ pieces on ${\mathcal G}'^{(\alpha)}$ with  weight $\tilde{Y}$ 
for pieces on odd sites $1,3,5,\cdots, 2\alpha-3$, $\tilde{P}$ for 
pieces on even sites $2,4,\cdots, 2\alpha - 4$ and the special weights $P/(Y(Y+P))$
for pieces on the site $2\alpha-2$, $1/(Y+P)$ for pieces on the site $2\alpha-1$ and $0$ for pieces on the site $2\alpha$.
With these definitions (and $k'^{(\alpha)}_{0}\equiv 1$), we have
\begin{equation*}
\sum_{m=0}^\alpha (-1)^m x'^{(\alpha)}_m\,  k'^{(\alpha)}_{n-m}=0\ \hbox{for all integers}\ n
\end{equation*} 
where $x'^{(\alpha)}_m$ enumerates configurations of $m$ hard pieces on ${\mathcal G}'^{(\alpha)}$
with weight $Y$ (resp.\ $P$) per piece sitting on an odd (resp.\ even) site (in particular $x'^{(\alpha)}_0=1$).
Setting $\alpha=i-1$,
the above equation becomes
\begin{equation*}
\sum_{m=0}^{i-1} (-1)^m x'^{(i-1)}_m\,  k'^{(i-1)}_{n-m}=0\ \hbox{for all integers}\ n  .
\end{equation*}
Again, from their heap interpretation, it is clear that $k_n=k_n'^{(i-1)}$ for $n=0,1,\cdots, i$ as well as for 
$n=-1,-2,\cdots,-(i-3)$. For $2\leq n\leq i$, all the $k'^{(i-1)}_p$'s appearing in the above formula may
thus be
replaced by $k_p$'s and \eqref{eq:linearC} follows. For $n=1$, the only term for which this substitution fails
is for $m=i-1$, since $k'^{(i-1)}_{-(i-2)}=k_{-(i-2)}-\tilde{P}^{i-3}(\tilde{P}-P/(Y(Y+P)))$ (the two indeed differ
by the contribution of the last piece, at position $2(i-2)$, in the the heap made of one piece on each even site from $2$ to $2(i-2)$). This
explains the right hand side $(-1)^{i-1}\tilde{P}^{i-3}(\tilde{P}-P/(Y(Y+P))) x'^{(i-1)}_{i-1}=(-1)^{i-1}P^{i-1}/Y^{i-2}$
 in the first line of \eqref{eq:linearD} (note that $x'^{(i-1)}_{i-1}=Y^{i-2}(Y+P)$).
For $n=i+1$, the only term for which this substitution fails is for $m=0$, since $k'^{(i-1)}_{i+1}=k_{i+1} - Y P^{i-1}(Y+P)$ (the two indeed differ
by the contribution of the two heaps made of one piece on site $1$, one piece on each even site from $2$ to $2(i-1)$
and a last piece at position $2i-1$ or $2i$). This
explains the right hand side $Y P^{i-1}(Y+P)$ 
in the second line of \eqref{eq:linearD}.

\section*{Acknowledgements} We warmly thank P. Di Francesco and R. Kedem for useful discussions. The work of
\'EF was partly supported by the ANR grant  
``Cartaplus'' 12-JS02-001-01 and the ANR grant ``EGOS'' 12-JS02-002-01.

\bibliographystyle{plain}
\bibliography{AB2p}

\end{document}